\documentclass{amsart}
\usepackage[margin=1in]{geometry} 
\usepackage{graphicx} 
\usepackage{amsmath}
\usepackage{numprint}
\usepackage{comment}
\usepackage{subcaption}
\usepackage{amsthm}
\usepackage{verbatim}
\usepackage{amssymb}
\usepackage{quiver}
\usepackage{xcolor}

\usepackage{float}
\usepackage[shortlabels]{enumitem}
\usepackage{quiver}
\usepackage{soul}
\usepackage[backend=biber,style=alphabetic]{biblatex}
\usepackage{tikz-cd}
\usetikzlibrary{nfold}
\usepackage[normalem]{ulem}
\usepackage{biblatex}
\usepackage{amsfonts}
\usepackage{hyperref}
\usepackage{adjustbox}
\renewbibmacro{in:}{}
\usepackage{cleveref}
\newtheorem{manualtheoreminner}{Theorem}
\newenvironment{manualtheorem}[1]{%
  \IfBlankTF{#1}
    {\renewcommand{\themanualtheoreminner}{\unskip}}
    {\renewcommand\themanualtheoreminner{#1}}%
  \manualtheoreminner
}{\endmanualtheoreminner}
\newtheorem{theorem}{Theorem}[section] 
\newtheorem{lemma}[theorem]{Lemma}

\newtheorem{question}[theorem]{Question}
\usepackage{mathtools}
\newtheorem{corollary}[theorem]{Corollary}

\newtheorem{remark}[theorem]{Remark}
\newtheorem{proposition}[theorem]{Proposition}

\newtheorem{example}[theorem]{Example}
\theoremstyle{definition}
\newtheorem{definition}[theorem]{Definition}

\newcommand{\il}[1]{\textit{#1}}
\newcommand{\g}[1]{\textbf{#1}}

\usepackage{xcolor}
\usepackage{algorithm}
\usepackage{algpseudocode}

\definecolor{editorcolor}{RGB}{255, 102, 0}

\title{Interpretations of the $\g{cd}$-index}

\begin{document}
\author[D. Guyer]{Dan Guyer}
\address[D. Guyer]{Department of Mathematics, University of Washington, Seattle, WA.}
\email{dguyer@uw.edu}

\author[A. Waugh]{Alex Waugh}
\address[A. Waugh]{Department of Mathematics, New Mexico State University, Las Cruces, NM.}
\email{awaugh@nmsu.edu}

\begin{abstract}
    We extend the recent work of Caster--Guyer--Samper (2026+) to establish a combinatorial interpretation of the \g{cd}-index for all $S$-partitionable posets, which includes all polytopes and $S$-shellable spheres. This interpretation is inductive on rank, but does not require passing to the order complex of the poset. We also provide topological interpretations of each \g{cd}-monomial. Furthermore, we establish a recurrence on the \g{cd}-indices of shelling components of the $d$-simplex that does not require the use of Andr\'e permutations, a unimodality result for \g{cd}-words ending in \g{d} for shelling components of the $d$-simplex and the $d$-cube, and provide two families of subdivisions of the boundary of the $d$-simplex whose sequence of \g{cd}-indices is unbounded on all but the monomial $\g{c}^{d+1}$ and all but the two monomials $\g{c}^{d+1}$ and $\g{d}\g{c}^{d-1}$.
\end{abstract}

\maketitle

\section{Introduction}

A central question in the study of geometric combinatorics is to classify which vectors in $\mathbb{R}^d$ arise as the face vectors (i.e. $f$-vectors) of $d$-dimensional polytopes. To better understand the combinatorial information encoded by the $f$-vector, we choose to study the finer flag $f$-vector and the closely-related invariant, the $\g{cd}$-index. The \textit{flag $f$-vector} of a $d$-polytope $P$ records the number of chains of faces of $P$. Namely, for each subset $\{i_1<i_2<\cdots<i_k\}\subseteq \{0,1,\ldots,d-1\}$
the flag $f$-vector records the number of chains $F_1 \subseteq F_2 \subseteq \cdots \subseteq F_k$,
where each $F_j$ is an $i_j$-dimensional face of $P$. Hence, the flag $f$-vector refines the usual $f$-vector with $d$ entries into a vector with $2^d$ entries.

Bayer and Billera~\cite{Billera_1983} showed that it is unnecessary to keep track of all $2^d$ entries of the flag $f$-vector of a $d$-polytope. Indeed, the Bayer--Billera relations imply that the affine dimension of the flag $f$-vectors of all Eulerian posets is $F_{d+1}-1$, where $F_{d+1}$ denotes the $(d+1)^{\text{th}}$ Fibonacci number. An encoding of the flag $f$-vector with only $F_{d+1}$ many entries was later introduced by Fine and further developed by Bayer and Klapper~\cite{Bayer1991}. This encoding, known as the \g{cd}-index, is a polynomial in noncommuting variables $\g{c}$ and $\g{d}$ of degree $1$ and $2$, respectively. Computational methods have been applied to existing inequalities on the flag $f$-vectors to derive new results on the $f$-vectors of polytopes in small dimensions by the work of Meisinger, Kleinschmidt and Kalai~\cite{Meisinger2000} and the very recent preprint of De Loera et~al.~\cite{deloerazheng2026}. Thus, to study the possible $f$-vectors of polytopes, we examine their flag $f$-vectors and \g{cd}-indices, which encode deeper combinatorial structure.

Understanding the coefficients of the \g{cd}-index remains challenging. One of the earliest works in this direction was that of Purtill~\cite{Purtill_1993}, who established the nonnegativity of the coefficients of the \g{cd}-index for certain classes of polytopes which include all simplicial and all simple polytopes. To do so, he used Andr\'e permutations of Foata--Sh\"utzenberger~\cite{Foata} to provide a combinatorial interpretation of the \g{cd}-index of a simplex. In this 1993 work, Purtill wrote 
\begin{quote}
\begin{flushleft}
``However, the \g{cd}-index appears to be a very subtle invariant; it is not even easy to calculate the \g{cd}-index of a simplex without some sort of general theory... Therefore, it would be nice to get some sort of combinatorial interpretation for the \g{cd}-index."
\end{flushleft}
\end{quote}

An algebraic interpretation of the \g{cd}-index was later established by Karu~\cite{Karu_2006} for all (regular) CW-spheres and a geometric interpretation was found by Lee~\cite{Lee2010} for all polytopes. However, a combinatorial interpretation has remained elusive. In fact, Margaret Bayer wrote the following in her 2021 survey~\cite{Bayer_2021}, 

\begin{quote}
\begin{flushleft}
    ``There are many specific open questions, but the overriding issue is to find combinatorial interpretations for the coefficients beyond those in special cases mentioned here."
\end{flushleft}
\end{quote}

In this article, we simplify the work of Caster--Guyer--Samper~\cite{caster2026} to provide a combinatorial interpretation of the $\g{cd}$-index for all $S$-partitionable posets (see Definition~\ref{def:S_part}). 
While their construction is combinatorial in nature, it does not provide an explicit answer to the question of ``What do the coefficients of the $\g{cd}$-index count?". For each $S$-partition $\nu$ of an Eulerian poset $\Omega$, we define a rooted tree $T_\nu(\Omega)$ which will be directed from the leaves to the root with some edges having a label from the set $\{1,2\}$ (see Figures~\ref{fig:Pentagon} and~\ref{fig:Cube}). This tree will be recursively defined and will provide an organizational tool for interpreting the recursion in Theorem~3.4 of Caster--Guyer--Samper~\cite{caster2026}. With this in mind, our main result can then be stated as follows.

\begin{manualtheorem}{A}
    \textit{If $\Omega$ is an $S$-partitionable poset, then the number of maximal monotone paths in $T_\nu(\Omega)$ with a weight sequence of $1$'s and $2$'s equals the coefficient of the $\g{cd}$-word of $\Omega$ obtained by replacing $1$ with $\g{c}$ and $2$ with $\g{d}$. } 
\end{manualtheorem}

As a consequence, we recover a new way to view Stanley's proof of the nonnegativity of the $\g{cd}$-index for all $S$-shellable spheres, and provide a way to visualize the formulas of Ehrenborg--Readdy~\cite{Ehrenborg1998} on the $\g{cd}$-index of the prism, pyramid and bipyramid of a polytope. Furthermore, we use these interpretations to establish a recurrence on the $\g{cd}$-indices of shelling components of the $d$-simplex that does not require the use of Andr\'e permutations.
We establish a unimodality result for \g{cd}-words ending in $\g{d}$ for shelling components of the $d$-simplex and the $d$-cube. And finally, we provide two families of subdivisions of the boundary of the $d$-simplex whose sequence of \g{cd}-indices is unbounded on all \g{cd}-monomials except $\g{c}^{d+1}$ and all \g{cd}-monomials except the two monomials $\g{c}^{d+1}$ and $\g{d}\g{c}^{d-1}$.

As a consequence of this interpretation of the $\g{cd}$-index, it becomes apparent that instead of looking at chains of faces directly, it suffices to keep track of the gluing data of the poset or CW complex when computing the $\g{cd}$-index. This topological perspective of the $\g{cd}$-index, in terms of gluing data, suggests that one may expect that the $\g{cd}$-index behaves well under certain topological operations. We build on Lemma~1.1 of Stanley~\cite{Stanley1994flag} by giving a topological interpretation of right multiplication by the monomial $\g{d}$. Furthermore, we provide topological models for every $\g{cd}$-monomial.

\section{Preliminaries}
Throughout, we work with finite graded posets, which we will often denote as $\Omega$, that have a unique minimum element $\hat{0}$ of rank $0$ and a unique maximum element $\hat{1}$ of rank $d+1$. If $\Omega$ is the face poset of a regular CW complex, we will typically not distinguish between the face poset and the CW complex itself. For further references to poset theory and the $\g{cd}$-index, we refer the reader to Sections 3.16 and 3.17 of Stanley's \textit{Enumerative Combinatorics} Volume 1~\cite{Stanley_ec1} and Bayer's recent survey~\cite{Bayer_2021}. For content related to polytopes and simplicial complexes, we refer the reader to Ziegler~\cite{Ziegler_1995}.

\begin{definition}
For a finite graded poset $\Omega$, the \textit{order complex} $\Delta(\Omega)$ is the simplicial complex whose faces are given by all chains, including the empty chain, of the elements of $\Omega\setminus \{\hat{0},\hat{1}\}$.
\end{definition}

 For a regular CW complex $\Omega$, the \textit{geometric realization of its order complex}, $||\Delta(\Omega)||$ is homeomorphic to the barycentric subdivision of $\Omega$, and hence to $\Omega$ itself~\cite[pg.~80]{Lundell_1969}. Thus, the order complex will recover the same topological information of $\Omega$, while turning the poset into a simplicial complex. For the purposes of this article, we will be primarily interested in studying the face posets of polytopes, regular CW spheres, and regular CW manifolds. 

\begin{definition}
A finite graded poset $\Omega$ with $\hat{0}$ and $\hat{1}$ is \textit{Eulerian} if every interval $[x,y]$ with $x<y$ has the same number of elements of odd rank as elements of even rank. Equivalently, a graded poset with $\hat{0}$ and $\hat{1}$ is \textit{Eulerian} if its Möbius function satisfies $\mu([x,y]) = (-1)^{\,r(y)-r(x)}$ for all intervals $[x,y]$.
\end{definition}

Thus, if a poset is Eulerian, then the Euler characteristic alone cannot distinguish this poset from the face poset of a sphere. In many cases, one would like to further impose the condition that the order complex of the poset in question is Cohen--Macaulay (see Bj\"orner~\cite{bjorner2014} for a survey), which guarantees that the reduced homology of $||\Delta(\Omega)||$ is concentrated in the top dimension. Hence all regular CW spheres are both Eulerian and Cohen--Macaulay; any poset that is both Eulerian and Cohen--Macaulay will be called \textit{Gorenstein*}. All Gorenstein* posets have a nonnegative \g{cd}-index by the main result of Karu~\cite{Karu_2006}.  

When thinking about this condition of being Eulerian geometrically and topologically, one has to be particularly careful about the elements $\hat{0},\hat{1}$ in the visualization. Namely, recall that the order complex of $\Omega$ does not include chains containing either $\hat{0}$ or $\hat{1}$. Therefore, visually, it is often convenient to not include either element in one's depiction. Thus, as in Lee~\cite{Lee2010}, when one considers the $\g{cd}$-index of a polytope, the interior of the polytope represents $\hat{1}$ in order to work with its spherical boundary. In some texts, as in Stanley~\cite{Stanley1994flag}, one may use $P_1(\Omega)$ to denote the poset $\Omega \cup \{\hat{1}\}$ for emphasis.

\begin{definition}
A finite graded poset $\Omega$ with $\hat{0}$ and $\hat{1}$ is \textit{semi-Eulerian} if every interval $[x,y]\neq [\hat{0},\hat{1}]$ with $x<y$ has the same number of elements of odd rank as elements of even rank. Equivalently, a graded poset with $\hat{0}$ and $\hat{1}$ is \textit{semi-Eulerian} if its Möbius function satisfies $\mu([x,y]) = (-1)^{\,r(y)-r(x)}$ for all intervals $[x,y]\neq [\hat{0},\hat{1}]$.
\end{definition}

Thus, semi-Eulerian posets are bounded graded posets whose Euler characteristic agrees with that of a sphere when taken on all proper intervals. Examples of semi-Eulerian posets include face posets of regular CW manifolds or, more generally, of homology manifolds. 

\begin{definition}
    We say that $\Omega$ is \textit{near-Eulerian} if it can be obtained from an Eulerian poset by deleting a single coatom.
\end{definition} 

Near-Eulerian posets arise naturally when building a recurrence for the $\g{cd}$-index, as done in Stanley~\cite{Stanley1994flag}. One such example is any $d$-dimensional CW ball, as one can glue a single top-dimensional cell along its boundary to make it into a $d$-dimensional sphere. Next, we introduce several poset-related definitions.

\begin{definition}
For $p\in \Omega$, we write $\overline{p} \coloneq [\hat{0},p]$ for the (down) \textit{closure} of $p$, and for a collection $\mathcal{C}\subseteq \Omega$, we define $
\overline{\mathcal{C}} \coloneq \bigcup_{p\in \mathcal{C}} \overline{p}.$
\end{definition}

\begin{definition}
    For a near-Eulerian poset $\Omega$, the closure of the union of all $y$ with $[y,\hat{1}]$ a three-element chain is called the \textit{boundary} of $\Omega$.
\end{definition}

Each such $y$ is a cocoatom that is contained in a single coatom. Topologically, this corresponds to all ridges that are contained in a single facet. For CW manifolds with boundary, the closure of these ridges exactly makes up the boundary. If $\Omega$ is near-Eulerian, then the associated Eulerian poset can be reconstructed by adjoining a new coatom $\tau$ that covers every $y$ with $[y,\hat{1}]$ being a three-element chain.

\begin{definition}
This operation of adjoining a new coatom $\tau$ that covers every $y$ with $[y,\hat{1}]$ a three-element chain is the \textit{semisuspension} of $\Omega$, and is denoted as $\tilde{\Sigma}\Omega$.
\end{definition}

The semisuspension operation will often be used to recover a sphere, or more generally an Eulerian poset, from something that is not. For example, we will often consider a collection of facets of a $(d+1)$-polytope that forms a $d$-ball on its boundary. A $d$-ball is not Eulerian. However, by gluing in a single $d$-cell along its boundary, one recovers a $d$-sphere, which is now Eulerian. This operation of gluing in a single $d$-cell, is precisely the semisuspension operation in this example.

Next, we wish to record all elements of $\Delta(\Omega)$ as a polynomial. This polynomial will be in noncommuting variables $\g{a}$ and $\g{b}$ of degree $1$. For a graded poset of $\Omega$ of rank $d+1$ with rank function $\rho$ and $K\subseteq [d]$, the $K$-flag number $f_K(\Omega)$ denotes the number of chains $\{x_1 \prec x_2 \prec \hdots \prec x_k:x_i\in \Omega\}$ whose \textit{rank set} is $\{\rho(x_i) : 1 \le i \le k\}=K$. Using $u_K$ to denote the word $u_K = u_1 \cdots u_d$ with $u_i = \g{a}$ if $i\notin K$ and $u_i = \g{b}$ if $i\in K$, we recover the \textit{flag $f$-polynomial of $\Omega$}, also known as the \textit{chain polynomial of $\Omega$}, as follows:
$$\Upsilon_\Omega(\g{a},\g{b}) = \sum_{K\subseteq [d]} f_K(\Omega)\, u_K.$$

However, to obtain more symmetry and structure among its coefficients, one often works with the \textit{flag $h$-vector} $(h_K)_{K\subseteq [d]}$, defined by the invertible transformation
$$ h_K(\Omega) = \sum_{T\subseteq K} (-1)^{|K\setminus T|} f_T(\Omega).$$

Again, wishing to store the entries of the flag $h$-vector as the coefficients of a polynomial, we obtain the \textit{flag $h$-polynomial of $\Omega$}, also known as the \textit{\g{ab}-index of $\Omega$}, as follows:
$$
\Psi_\Omega(\g{a},\g{b}) = \sum_{K\subseteq [d]} h_K(\Omega)\, u_K, 
\qquad
$$
The flag $f$- and flag $h$-polynomials are related by the (invertible) transformation
$$
\Upsilon_\Omega(\g{a}-\g{b},\g{b}) = \Psi_\Omega(\g{a},\g{b}).
$$

In the study of simplicial complexes, one makes a similar transformation between the $f$-polynomial and $h$-polynomial. Moreover, for Eulerian simplicial complexes, there are the famous Dehn--Sommerville equations which say that $h_i=h_{d-i}$ for all $0\leq i \leq d$. By the results of Bayer--Billera~\cite{Billera_1983,Billera_1985}, the $2^d$ entries of the flag $h$-polynomial also exhibit analogous symmetry. The \il{generalized Dehn--Sommerville equations}, introduced by Bayer-Billera~\cite{Billera_1983,Billera_1985}, describe the affine span of all flag $f$-vectors of Eulerian posets of rank $d+1$. 

\begin{theorem}[Bayer--Billera~\cite{Billera_1985}]\label{thm:Bayer-Billera}
    Let $\Omega$ be an Eulerian poset of rank $d+1$. Then $h_{K}=h_{[d]\setminus K}$ for every subset $K\subseteq [d]$.
\end{theorem}

We remark that these equations of $h_{K}=h_{[d]\setminus K}$ are different from those stated in Theorem~2.1 of Bayer--Billera~\cite{Billera_1985}, but are equivalent. Because of these relations, one can re-express the $2^d$ entries of the flag $h$-polynomial into a polynomial in noncommuting variables $\g{c}$ and $\g{d}$ of degrees $1$ and $2$ respectively with just $F_{d+1}\approx 1.618^d$ entries as shown by Fine and Bayer--Klapper~\cite{Bayer1991}.

\begin{definition}\label{defn:cdIndex}
Let $\Omega$ be a finite graded poset of rank $d+1$. The \emph{$\g{cd}$-index} of $\Omega$ is the unique polynomial $\Phi_\Omega(\g{c},\g{d})$ with $\g{c}=\g{a}+\g{b}$ and $\g{d}=\g{ab}+\g{ba}$ such that
\[
\Phi_\Omega(\g{a}+\g{b},\g{ab}+\g{ba}) = \Psi_\Omega(\g{a},\g{b}),
\]
whenever such a polynomial exists. 
\end{definition}
We may often abuse notation and write $\Phi(\Omega)$ as shorthand for $\Phi_\Omega(\g{c},\g{d})$. If $\Omega$ has a $\g{cd}$-index, then $\Phi_\Omega(\g{c},\g{d})$ is unique and homogeneous of degree $d$ in $\mathbb{Z}\langle \g{c},\g{d}\rangle$, assigning degree $1$ to $\g{c}$ and degree $2$ to \g{d}. The number of \g{cd}-monomials (i.e. $\g{cd}$-words) of degree $d$ equals the $(d+1)^{\text{th}}$ Fibonacci number $F_{d+1}$ defined with initial conditions $F_1=F_2=1$. Furthermore, the following theorem of Bayer and Klapper~\cite{Bayer1991}, shows that a $\g{cd}$-index exists exactly when the generalized Dehn--Sommerville equations hold. 
\begin{theorem}[Bayer--Klapper~\cite{Bayer1991}]\label{thm: BayerKlapper}

A graded poset has a $\g{cd}$-index if and only if its flag $f$-vector satisfies the generalized Dehn--Sommerville equations. In particular, every Eulerian poset has a $\g{cd}$-index.
\end{theorem}

\section{A Combinatorial Interpretation}

In this section, we will provide a combinatorial interpretation of the $\g{cd}$-index for all $S$-partitionable posets. This construction illustrates that one need not pass to the order complex of an $S$-partitionable poset to compute its $\g{cd}$-index. In the specific cases of polytopes, and more generally $S$-shellable spheres, we show that the $\g{cd}$-index can be directly recovered from any $S$-shelling and thus any Brugesser--Mani line shelling~\cite{Bruggesser}. Our interpretation provides a concrete answer to the question of ``What do the coefficients of the $\g{cd}$-index count?" for $S$-partitionable posets. It is the number of maximal paths with a given sequence of weights in a weighted tree that records the $S$-partition. Each weight will be either $1$ or $2$ and each maximal path will contribute a single $\g{cd}$-word associated to this sequence of weights. A similar combinatorial interpretation readily follows for $SE$-partitionable posets, as defined in Caster--Guyer--Samper~\cite{caster2026}, but we expect such a description to be out of reach for more general Eulerian posets. This is because for simplicial complexes, at the present time, it is only known that the $h$-vector carries combinatorial meaning when the simplicial complex in question is partitionable; $S$-partitionability is the proposed analog of partitionability for Eulerian posets and the $\g{cd}$-index. To illustrate this new combinatorial description, we will also provide a few examples in this section.

\subsection{Our Combinatorial Interpretation}
\begin{figure}[h!]
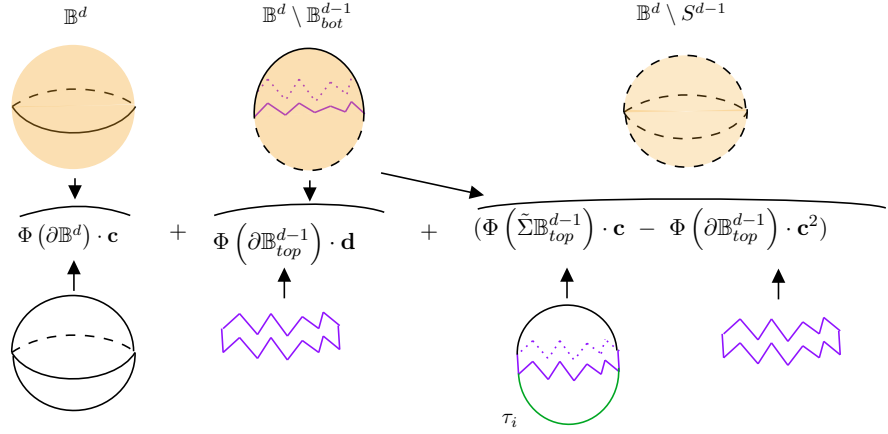

    \centering
    
\adjustbox{scale=0.8}{

\tikzset{every picture/.style={line width=0.75pt}} 



}

    \caption{A Motivating Heuristic for $S$-partitionability}
    \label{fig:heuristic}
\end{figure}

We begin by recalling the definition of $S$-partitionable posets, as defined in Caster--Guyer--Samper~\cite{caster2026}. As motivation, for $\Omega$ a regular CW $d$-sphere, one should think of an $S$-partition of $\Omega$ as a partition of its cells consisting of a single $d$-ball (the initial class), a single $d$-ball with all of its boundary removed (the terminal class), and any number of $d$-balls each with a $(d-1)$-ball removed from its boundary (the ordinary classes). Inductively, the $\g{cd}$-index of $\Omega$ will be computed from the boundary of the initial class as well as from each ordinary class in two steps. After removing the single top-dimensional cell from the ordinary class, one takes the closure to get a $(d-1)$-ball $\Gamma_\sigma$. Then, there are two spheres that can be associated to $\Gamma_\sigma$: the $(d-2)$-sphere $\partial\Gamma_{\sigma}$ and the $(d-1)$-sphere $\tilde{\Sigma}\Gamma_\sigma$. We illustrate how the \g{cd}-index can be computed from each $S$-partition class in Figure~\ref{fig:heuristic}

Then, the inductive assumptions in the definition of $S$-partitionable posets allow the spaces to be further broken down. By the results of Caster--Guyer--Samper~\cite{caster2026}, any line shelling of a polytope, any $S$-shelling of a sphere, or partitioning of simplicial sphere provides such a decomposition. In its fullest generality, the definition of $S$-partitionability is as follows:

\begin{definition}\label{def:S_part}
Let $\Omega$ be an Eulerian poset of rank $d+1$. We say that $\Omega$ is \textit{\il{S}-partitionable} if either $\Omega=\{\emptyset\}$ or $\Omega$ has rank at least 2 and there exists a partition of the elements of $\Omega$ with a partition class $P_\sigma \subseteq \overline{\sigma}$ for each coatom $\sigma$ of $\Omega$. Furthermore, for $\Gamma_\sigma \coloneq \overline{P_\sigma \setminus \{\sigma\}}$ we require:
        \begin{enumerate}[(a)]
            \item There is a unique \textit{initial} coatom $\sigma_1$ for which $P_{\sigma_1}=\overline{\sigma_1}$ and we have that $P_{\sigma_1}\setminus \{\sigma_1\}=\partial \overline{\sigma_1}$ is \il{S}-partitionable of rank $d$.
            \item There is a unique \textit{terminal} coatom $\sigma_t$ for which $P_{\sigma_t}= \{\sigma_t\}$ is a one element set. 
            \item For every other coatom $\sigma$, $\Gamma_\sigma$ is near-Eulerian of rank $d$ with $\Gamma_\sigma = (P_\sigma \setminus \{\sigma\}) \uplus \partial \Gamma_\sigma$, and $\tilde{\Sigma}\Gamma_\sigma$ is \il{S}-partitionable of rank $d$ such that the initial coatom is the added coatom $\tau_\sigma$.
        \end{enumerate}
\end{definition}
We will say that the coatoms and partition classes described in part (c) of Definition~\ref{def:S_part} are \textit{ordinary}. 
\begin{remark}\label{rem:boundary_is_S_part}
    For every ordinary coatom $\sigma$ of an \il{S}-partition, the poset $\partial \Gamma_\sigma$ is \il{S}-partitionable of rank $d-1$ because of property $(a)$ of Definition~\ref{def:S_part} applied to $\tilde{\Sigma}\Gamma_\sigma$.
\end{remark}

\begin{remark}\label{rem:facets_part}
    Observe that because $P_\sigma \subseteq \overline{\sigma}$, we find that $\sigma \in P_\sigma$ for each coatom $\sigma$ of $\Omega$. This follows from the fact that $\sigma \in \overline{\sigma'}$ if and only if $\sigma=\sigma'$ for two coatoms $\sigma,\sigma'$. 
\end{remark} 

In Caster--Guyer--Samper~\cite{caster2026}, they proved a recursive formula for the $\g{cd}$-index of $\Omega$ based on the $\g{cd}$-index of each type of coatom. Let $\mathcal{C}_{\Omega}(\sigma)$ denote the contribution of $\sigma$ to the \textbf{cd}-index of $\Omega$. Then this recursive formula is given as follows: 

$$\Phi(\Omega)=\Phi(\partial \overline{\sigma_1})\cdot \g{c} + \sum_{\substack{\sigma \text{ ordinary}\\ \text{in $\Omega$} }}\Phi(\partial \Gamma_\sigma)\cdot \mathbf{d} + \sum_{\substack{\sigma \text{ ordinary}\\ \text{in $\Omega$} }} \ \sum_{\substack{\omega \text{ ordinary}\\ \text{in $\tilde{\Sigma}\Gamma_\sigma$} }}\mathcal{C}_{\tilde{\Sigma}\Gamma_\sigma}(\omega) \cdot \g{c}.$$

 For the first term, we know that $\partial\overline{\sigma_1}$ is $S$-partitionable of rank $d$, while $\Omega$ has rank $d+1$. Thus, by induction, $\partial \overline{\sigma}_1$ has a $\g{cd}$-index that can be recursively computed; take this $\g{cd}$-index and multiply it on the right by $\g{c}$. Next, we move onto the ordinary coatoms. By Remark~\ref{rem:boundary_is_S_part}, $\partial\Gamma_\sigma$ is $S$-partitionable of rank $d-1$. Thus, by induction, $\partial \Gamma_\sigma$ has a $\g{cd}$-index that can be recursively computed; take this $\g{cd}$-index and multiply it on the right by $\g{d}$. For the final term, we first recognize that $\tilde{\Sigma}\Gamma_\sigma$ is $S$-partitionable of rank $d$. However, we had to ``artificially" add in $\tau_\sigma$ to make it so. Therefore, we will want to consider the contributions from all coatoms except the initial one in the semisuspension of $\Gamma_\sigma$; by induction these can be recursively computed. To establish a combinatorial interpretation of the $\g{cd}$-index, we will exhibit a recursively defined weighted directed graph that is a tree, and show that the coefficients of the $\g{cd}$-index are in bijection with the maximal directed paths of this tree, where each $\g{cd}$-word will correspond to a different weight sequence of a maximal directed path.  Define $T_\nu(\Omega)$, the \textit{Tree of the $S$-partition $\nu$ of $\Omega$}, recursively from the following rules.

Start with a node representing $\Omega$, the root of $T_{\nu}(\Omega)$. Then, for each non-terminal $S$-partition class of $\Omega$, draw a new node with an edge directed towards $\Omega$ to indicate that the $\g{cd}$-index of $\Omega$ depends on each of these contributions; these edges will break down $\Omega$ into its constituent $S$-partition classes and will not have a weight. For each such $S$-partition class of $\Omega$, draw additional nodes and edges as follows:
\begin{itemize}
    \item For the initial class of $\Omega$, draw a single node representing $\partial \overline{\sigma_1}$ with an edge directed from $\partial \overline{\sigma_1}$ to $\overline{\sigma_1}$ with weight $1$.
    \item For every ordinary class of $\Omega$, draw a node representing $\partial \Gamma_\sigma$ with an edge directed from $\partial \Gamma_\sigma$ to the $S$-partition class $P_\sigma$ of weight $2$, and a node for every ordinary $S$-partition class of $\tilde{\Sigma}\Gamma_\sigma$ with an edge directed towards $P_\sigma$ of weight $1$. 
\end{itemize}
 Repeat until every leaf of the tree represents $\{\emptyset\}$. Notice that all edges are directed towards the root, and the weights on every edge correspond to the difference of the ranks of the two sub-posets in question.

\textbf{Examples:} We include a couple of figures to illustrate this interpretation before giving its proof. 
\begin{figure}[h!]
    \centering
    
\adjustbox{scale = 0.95}{

\tikzset{every picture/.style={line width=0.75pt}} 

\begin{tikzpicture}[x=0.75pt,y=0.75pt,yscale=-1,xscale=1]

\draw  [color={rgb, 255:red, 245; green, 166; blue, 35 }  ,draw opacity=1 ][fill={rgb, 255:red, 255; green, 255; blue, 255 }  ,fill opacity=1 ] (115.46,77.99) -- (106.32,105.92) -- (76.93,105.86) -- (67.91,77.9) -- (91.71,60.67) -- cycle ;
\draw [color={rgb, 255:red, 245; green, 166; blue, 35 }  ,draw opacity=1 ][line width=1.5]    (140.33,80.33) -- (170.67,80.33) ;
\draw  [color={rgb, 255:red, 245; green, 166; blue, 35 }  ,draw opacity=1 ][fill={rgb, 255:red, 245; green, 166; blue, 35 }  ,fill opacity=1 ] (89.21,60.67) .. controls (89.21,59.29) and (90.33,58.17) .. (91.71,58.17) .. controls (93.1,58.17) and (94.21,59.29) .. (94.21,60.67) .. controls (94.21,62.05) and (93.1,63.17) .. (91.71,63.17) .. controls (90.33,63.17) and (89.21,62.05) .. (89.21,60.67) -- cycle ;
\draw  [color={rgb, 255:red, 245; green, 166; blue, 35 }  ,draw opacity=1 ][fill={rgb, 255:red, 245; green, 166; blue, 35 }  ,fill opacity=1 ] (112.96,77.99) .. controls (112.96,76.61) and (114.08,75.49) .. (115.46,75.49) .. controls (116.84,75.49) and (117.96,76.61) .. (117.96,77.99) .. controls (117.96,79.37) and (116.84,80.49) .. (115.46,80.49) .. controls (114.08,80.49) and (112.96,79.37) .. (112.96,77.99) -- cycle ;
\draw  [color={rgb, 255:red, 245; green, 166; blue, 35 }  ,draw opacity=1 ][fill={rgb, 255:red, 245; green, 166; blue, 35 }  ,fill opacity=1 ] (103.82,105.92) .. controls (103.82,104.54) and (104.94,103.42) .. (106.32,103.42) .. controls (107.7,103.42) and (108.82,104.54) .. (108.82,105.92) .. controls (108.82,107.3) and (107.7,108.42) .. (106.32,108.42) .. controls (104.94,108.42) and (103.82,107.3) .. (103.82,105.92) -- cycle ;
\draw  [color={rgb, 255:red, 245; green, 166; blue, 35 }  ,draw opacity=1 ][fill={rgb, 255:red, 245; green, 166; blue, 35 }  ,fill opacity=1 ] (74.43,105.86) .. controls (74.43,104.48) and (75.55,103.36) .. (76.93,103.36) .. controls (78.31,103.36) and (79.43,104.48) .. (79.43,105.86) .. controls (79.43,107.24) and (78.31,108.36) .. (76.93,108.36) .. controls (75.55,108.36) and (74.43,107.24) .. (74.43,105.86) -- cycle ;
\draw  [color={rgb, 255:red, 245; green, 166; blue, 35 }  ,draw opacity=1 ][fill={rgb, 255:red, 245; green, 166; blue, 35 }  ,fill opacity=1 ] (65.41,77.9) .. controls (65.41,76.51) and (66.52,75.4) .. (67.91,75.4) .. controls (69.29,75.4) and (70.41,76.51) .. (70.41,77.9) .. controls (70.41,79.28) and (69.29,80.4) .. (67.91,80.4) .. controls (66.52,80.4) and (65.41,79.28) .. (65.41,77.9) -- cycle ;
\draw  [color={rgb, 255:red, 245; green, 166; blue, 35 }  ,draw opacity=1 ][fill={rgb, 255:red, 245; green, 166; blue, 35 }  ,fill opacity=1 ] (137.83,80.33) .. controls (137.83,78.95) and (138.95,77.83) .. (140.33,77.83) .. controls (141.71,77.83) and (142.83,78.95) .. (142.83,80.33) .. controls (142.83,81.71) and (141.71,82.83) .. (140.33,82.83) .. controls (138.95,82.83) and (137.83,81.71) .. (137.83,80.33) -- cycle ;
\draw  [color={rgb, 255:red, 245; green, 166; blue, 35 }  ,draw opacity=1 ][fill={rgb, 255:red, 245; green, 166; blue, 35 }  ,fill opacity=1 ] (168.17,80.33) .. controls (168.17,78.95) and (169.29,77.83) .. (170.67,77.83) .. controls (172.05,77.83) and (173.17,78.95) .. (173.17,80.33) .. controls (173.17,81.71) and (172.05,82.83) .. (170.67,82.83) .. controls (169.29,82.83) and (168.17,81.71) .. (168.17,80.33) -- cycle ;
\draw  [color={rgb, 255:red, 245; green, 166; blue, 35 }  ,draw opacity=1 ][fill={rgb, 255:red, 245; green, 166; blue, 35 }  ,fill opacity=1 ] (138.17,130.33) .. controls (138.17,128.95) and (139.29,127.83) .. (140.67,127.83) .. controls (142.05,127.83) and (143.17,128.95) .. (143.17,130.33) .. controls (143.17,131.71) and (142.05,132.83) .. (140.67,132.83) .. controls (139.29,132.83) and (138.17,131.71) .. (138.17,130.33) -- cycle ;
\draw  [color={rgb, 255:red, 245; green, 166; blue, 35 }  ,draw opacity=1 ][fill={rgb, 255:red, 245; green, 166; blue, 35 }  ,fill opacity=1 ] (167.83,130.67) .. controls (167.83,129.29) and (168.95,128.17) .. (170.33,128.17) .. controls (171.71,128.17) and (172.83,129.29) .. (172.83,130.67) .. controls (172.83,132.05) and (171.71,133.17) .. (170.33,133.17) .. controls (168.95,133.17) and (167.83,132.05) .. (167.83,130.67) -- cycle ;
\draw    (155.5,124.67) -- (155.2,93.67) ;
\draw [shift={(155.17,90.67)}, rotate = 89.44] [fill={rgb, 255:red, 0; green, 0; blue, 0 }  ][line width=0.08]  [draw opacity=0] (8.93,-4.29) -- (0,0) -- (8.93,4.29) -- cycle    ;
\draw  [dash pattern={on 4.5pt off 4.5pt}]  (172.97,118.3) .. controls (185.51,105.55) and (202.55,106.41) .. (212,122) ;
\draw [shift={(170.83,120.67)}, rotate = 309.87] [fill={rgb, 255:red, 0; green, 0; blue, 0 }  ][line width=0.08]  [draw opacity=0] (8.93,-4.29) -- (0,0) -- (8.93,4.29) -- cycle    ;
\draw  [color={rgb, 255:red, 245; green, 166; blue, 35 }  ,draw opacity=1 ][fill={rgb, 255:red, 245; green, 166; blue, 35 }  ,fill opacity=1 ] (207.33,130.83) .. controls (207.33,129.45) and (208.45,128.33) .. (209.83,128.33) .. controls (211.21,128.33) and (212.33,129.45) .. (212.33,130.83) .. controls (212.33,132.21) and (211.21,133.33) .. (209.83,133.33) .. controls (208.45,133.33) and (207.33,132.21) .. (207.33,130.83) -- cycle ;
\draw    (211,172) -- (210.7,141) ;
\draw [shift={(210.67,138)}, rotate = 89.44] [fill={rgb, 255:red, 0; green, 0; blue, 0 }  ][line width=0.08]  [draw opacity=0] (8.93,-4.29) -- (0,0) -- (8.93,4.29) -- cycle    ;
\draw [color={rgb, 255:red, 245; green, 166; blue, 35 }  ,draw opacity=1 ][line width=1.5]    (299.5,80) -- (329.83,80) ;
\draw  [color={rgb, 255:red, 245; green, 166; blue, 35 }  ,draw opacity=1 ][fill={rgb, 255:red, 245; green, 166; blue, 35 }  ,fill opacity=1 ] (297,80) .. controls (297,78.62) and (298.12,77.5) .. (299.5,77.5) .. controls (300.88,77.5) and (302,78.62) .. (302,80) .. controls (302,81.38) and (300.88,82.5) .. (299.5,82.5) .. controls (298.12,82.5) and (297,81.38) .. (297,80) -- cycle ;
\draw  [color={rgb, 255:red, 245; green, 166; blue, 35 }  ,draw opacity=1 ][fill={rgb, 255:red, 255; green, 255; blue, 255 }  ,fill opacity=1 ] (327.33,80) .. controls (327.33,78.62) and (328.45,77.5) .. (329.83,77.5) .. controls (331.21,77.5) and (332.33,78.62) .. (332.33,80) .. controls (332.33,81.38) and (331.21,82.5) .. (329.83,82.5) .. controls (328.45,82.5) and (327.33,81.38) .. (327.33,80) -- cycle ;
\draw    (314.33,179.33) -- (314.33,93.33) ;
\draw [shift={(314.33,90.33)}, rotate = 90] [fill={rgb, 255:red, 0; green, 0; blue, 0 }  ][line width=0.08]  [draw opacity=0] (8.93,-4.29) -- (0,0) -- (8.93,4.29) -- cycle    ;
\draw    (289.67,69.67) .. controls (275.67,89.5) and (276.33,182) .. (292.33,201) ;
\draw    (339.33,70.67) .. controls (358.67,89.67) and (359.67,171.33) .. (340,200.23) ;
\draw  [dash pattern={on 4.5pt off 4.5pt}]  (106.56,58.38) .. controls (125.27,47.79) and (155.71,56.37) .. (156.67,71.33) ;
\draw [shift={(104,60)}, rotate = 324.66] [fill={rgb, 255:red, 0; green, 0; blue, 0 }  ][line width=0.08]  [draw opacity=0] (8.93,-4.29) -- (0,0) -- (8.93,4.29) -- cycle    ;
\draw  [dash pattern={on 4.5pt off 4.5pt}]  (105.29,48.55) .. controls (134.86,37.81) and (301.71,47.63) .. (302.67,62.67) ;
\draw [shift={(102.33,50)}, rotate = 324.66] [fill={rgb, 255:red, 0; green, 0; blue, 0 }  ][line width=0.08]  [draw opacity=0] (8.93,-4.29) -- (0,0) -- (8.93,4.29) -- cycle    ;
\draw  [dash pattern={on 4.5pt off 4.5pt}]  (101.97,40.21) .. controls (131.79,29.29) and (301.71,36.96) .. (302.67,52) ;
\draw [shift={(99,41.67)}, rotate = 324.66] [fill={rgb, 255:red, 0; green, 0; blue, 0 }  ][line width=0.08]  [draw opacity=0] (8.93,-4.29) -- (0,0) -- (8.93,4.29) -- cycle    ;
\draw  [dash pattern={on 4.5pt off 4.5pt}]  (99.97,30.21) .. controls (129.79,19.29) and (299.71,26.96) .. (300.67,42) ;
\draw [shift={(97,31.67)}, rotate = 324.66] [fill={rgb, 255:red, 0; green, 0; blue, 0 }  ][line width=0.08]  [draw opacity=0] (8.93,-4.29) -- (0,0) -- (8.93,4.29) -- cycle    ;

\draw (140.17,100.9) node [anchor=north west][inner sep=0.75pt]  [font=\small]  {$1$};
\draw (192.67,151.9) node [anchor=north west][inner sep=0.75pt]  [font=\small]  {$1$};
\draw (204.67,179.23) node [anchor=north west][inner sep=0.75pt]  [font=\small,color={rgb, 255:red, 245; green, 166; blue, 35 }  ,opacity=1 ]  {$\{\emptyset\} $};
\draw (325,131.57) node [anchor=north west][inner sep=0.75pt]  [font=\small]  {$2$};
\draw (308.17,181.3) node [anchor=north west][inner sep=0.75pt]  [font=\small,color={rgb, 255:red, 245; green, 166; blue, 35 }  ,opacity=1 ]  {$\{\emptyset\} $};
\draw (361,129.9) node [anchor=north west][inner sep=0.75pt]  [font=\small]  {$\times 3$};

\end{tikzpicture}
}

    \caption{$T_\nu(\Omega)$ for the pentagon.}
    \label{fig:Pentagon}
\end{figure}

In Figure~\ref{fig:Pentagon}, we find that the pentagon has a $\g{cd}$-index of $\g{c}^2+3\g{d}$. This is because there is one path with weight sequence $1-1$ and three paths with weight sequence $2$. 

\begin{figure}[h!]
\centering
     
\adjustbox{scale = 0.95}{

\tikzset{every picture/.style={line width=0.75pt}} 



}
\caption{$T_\nu(\Omega)$ for the $3$-cube.}
    \label{fig:Cube}
\end{figure}

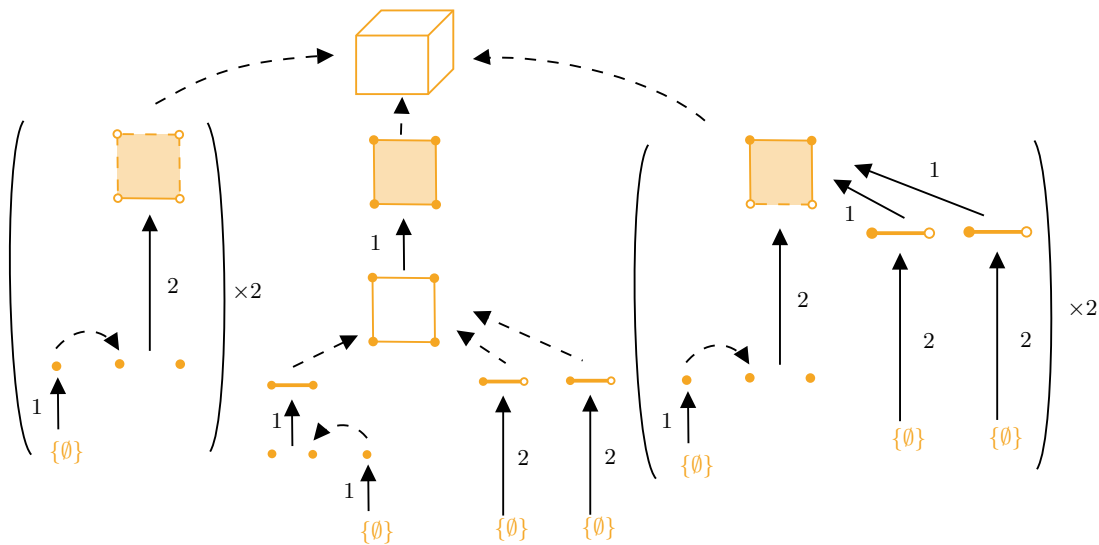
In Figure~\ref{fig:Cube}, we find that the $3$-cube has a $\g{cd}$-index of $\g{c}^3+4\g{cd}+6\g{dc}$. This is because there is one path with weight sequence $1-1-1$, four paths of weight sequence $1-2$ and six paths of weight sequence $2-1$. 

The construction of $T_\nu(\Omega)$ is built to organize and understand the recursion of Theorem~3.4 of Caster--Guyer--Samper~\cite{caster2026}. Therefore, the proof of Theorem~\ref{thm:combint} will verify that the recursion of Caster--Guyer--Samper~\cite{caster2026} is indeed recorded by $T_\nu(\Omega)$. 

\begin{theorem}\label{thm:combint}
    If $\Omega$ is an $S$-partitionable poset, then the number of maximal monotone paths in $T_\nu(\Omega)$ with a weight sequence of $1$'s and $2$'s equals the coefficient of the $\g{cd}$-word of $\Omega$ obtained by replacing $1$ with $\g{c}$ and $2$ with $\g{d}$. 
\end{theorem}

\begin{proof} We will proceed via induction on the rank of $\Omega$. The two base cases are as follows: the poset $\{\emptyset\}$, which has a $\g{cd}$-index of $1$ and the diamond poset which has two elements of rank one. The diamond poset can be modeled as $S^0$, which admits a single initial coatom and a single terminal coatom. In the definition of $T_{\nu}(\Omega)$, the terminal $S$-partition class is never included which yields the following for $T_\nu(S^0)$:
    
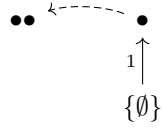
\begin{center}
\begin{tikzcd}
	{\bullet \bullet} & \bullet \\
	& {\{\emptyset\}}
	\arrow[curve={height=6pt}, dashed, from=1-2, to=1-1]
	\arrow["1", from=2-2, to=1-2]
\end{tikzcd}
\end{center}
From this, we compute a $\g{cd}$-index of $\g{c}$ as there is a single path of weight $1$. This indeed matches the \g{cd}-index of $S^0$, which one can compute directly from the flag $f$-polynomial.
Assume for the sake of induction that the result holds for all posets of rank at most $d-1$ and that $\Omega$ has rank $d$. Then, the recurrence from Caster--Guyer--Samper~\cite[Theorem~3.4]{caster2026} shows that the $\g{cd}$-index of $\Omega$ can be recursively computed via the sums of the contributions from each $S$-partition class:
$$\Phi(\Omega)=\Phi(\partial \overline{\sigma_1})\cdot \g{c} + \sum_{\substack{\sigma \text{ ordinary}\\ \text{in $\Omega$} }}\Phi(\partial \Gamma_\sigma)\cdot \mathbf{d} + \sum_{\substack{\sigma \text{ ordinary}\\ \text{in $\Omega$} }} \ \sum_{\substack{\omega \text{ ordinary}\\ \text{in $\tilde{\Sigma}\Gamma_\sigma$} }}\mathcal{C}_{\tilde{\Sigma}\Gamma_\sigma}(\omega) \cdot \g{c}.$$

Since $P_1$, the initial $S$-partition class, has a contribution of $\partial \overline{\sigma_1}\cdot \g{c}$, and $P_1$ is $S$-partitionable of rank $d-1$, its contribution is correct in $T_\nu(\Omega)$ via induction. If $\sigma$ is an ordinary coatom, then its directed edge with weight $2$ going from $\partial \Gamma_\sigma$ to its $S$-partition class also contributes the correct contribution by induction since $\partial \Gamma_\sigma$ is $S$-partitionable of rank $d-2$. Moreover, $\tilde{\Sigma}\Gamma_\sigma$ is $S$-partitionable of rank $d-1$. Therefore, by induction, one can exhibit the decomposition for $\tilde{\Sigma}\Gamma_\sigma$ and remove the subtree associated to the initial partition class $\tau_\sigma$ to obtain the correct contributions in the recursion, matching the construction of $T_\nu(\Omega)$. Finally, we remark that in the Caster--Guyer--Samper recurrence the terminal $S$-partition class contributes nothing to the sum, which is reflected in not including it in $T_\nu(\Omega)$. 
\end{proof}

\begin{remark}
    An analogous construction can be used to record the recursion of Theorem~4.2 of Caster--Guyer--Samper~\cite{caster2026} for $SE$-partitionable posets. We omit the proof for brevity. 
\end{remark}

From Theorem~\ref{thm:combint}, we immediately pick up a weakening of Theorem~3.1 of Stanley~\cite{Stanley1994flag} regarding the $\g{cd}$-indices of simplicial complexes. We shall use $\hat{\Phi}_k^d$ to denote the $\g{cd}$-index of the shelling component of the boundary of the $d$-simplex with a restriction face of size $k$.  

\begin{corollary}
    If $\Delta$ is a $d$-dimensional Eulerian partitionable simplicial complex, then $\Phi(\Delta)=\sum_{i=0}^dh_i(\Delta)\hat{\Phi}_i^d$.
\end{corollary}
\begin{proof}
    From Theorem~\ref{thm:combint}, the $\g{cd}$-index is determined from the sum of the contributions of each partition class since Eulerian partitionable simplicial complexes are $S$-partitionable by Corollary~3.11 of Caster--Guyer--Samper~\cite{caster2026}. There are finitely many possible partition classes of the boundary of the $d$-simplex, each one's contribution is denoted by $\hat{\Phi}_i^d$ where $i$ is the size of the restriction face. Since there $d+1$ options for $i$ and $h_i$ counts the number partition classes of each type, the proof follows. 
\end{proof}

By Theorem~3.8 of Ehrenborg--Hetyei~\cite{EhrenborgHetyei2000}, a total order of the facets of a cubical sphere is an $S$-shelling if and only if it is a $C$-shelling. Therefore, there is no need to specify the type of shelling for cubical spheres. We say that a shelling component of the $d$-cube is of type $(r,s)$ if there are exactly $r$ facets present without their antipodes, and $s$ pairs of antipodal facets in which both facets are present. For a shellable cubical complex, one counts the number of shelling components of each type within the $c$-vector ($\ldots,c_{r,s},\ldots$). The following result also immediately follows from Theorem~\ref{thm:combint}.

\begin{corollary}[Ehrenborg--Hetyei~\cite{EhrenborgHetyei2000}]
Let $\mathcal{C}$ be an $(n-1)$-dimensional shellable cubical sphere which has a shelling with $c$-vector $(\ldots,c_{r,s},\ldots,)$. Then the $\g{cd}$-index of $P_1(\mathcal{C})$ is given by 
$$\Phi(P_1(\mathcal{C}))=V_{n+1,1,0}+\sum_{\substack{r,s\\ r>0}}c_{r,s}(V_{n+1,r+1,s}-V_{n+1,r,s}).$$
\end{corollary}

\begin{proof}
    The quantity $V_{n+1,r+1,s}-V_{n+1,r,s}$ represents the $\g{cd}$-index contribution from each shelling component while $c_{r,s}$ counts the number of shelling components of each type. Since there is always a single initial partition class, this is counted by $V_{n+1,1,0}$.
\end{proof}

\subsection{Nonnegativity of Stanley's $S$-shellable spheres}
Via Theorem~\ref{thm:combint}, we recover a quick proof of Theorem~2.2 of Stanley~\cite{Stanley1994flag} on the nonnegativity of the $\g{cd}$-index for all $S$-shellable spheres. In Caster--Guyer--Samper~\cite{caster2026}, they verified that the inductive properties of $S$-shellability imply the inductive properties of $S$-partitionability. However, this observation alone does not shed light onto the steps of Stanley's proof. His proof starts by defining a sequence of CW spheres $\Omega_1,\Omega_2,\ldots,\Omega_{r-1}$. Then, he shows that the $\g{cd}$-index of $\Omega_1$ is nonnegative and that the $\g{cd}$-index of $\Omega_i$ is (coefficient-wise) less than or equal to the $\g{cd}$-index of $\Omega_j$ for $i<j$. In this subsection, we provide a proof of Stanley's result that quickly illustrates why both of these two steps hold. We first recall Stanley's definition of $S$-shellable spheres.

The motivation will be similar to that of $S$-partitionable posets when $\Omega$ is a $d$-sphere. One should think of an $S$-shelling as a total order on the facets that starts with a single $d$-ball (with all of its boundary), then adding some number of $d$-balls each with a $(d-1)$-ball removed from its boundary to indicate where it is glued to the rest of the complex, before finishing by adding on a single $d$-ball with all of its boundary removed. To continue the induction, it is required that the facets on the boundary of each added $d$-ball come equipped with their own $S$-shelling order. The definition of $S$-shellability is as follows:

\begin{definition}[Stanley~\cite{Stanley1994flag}]\label{defn:Sshell} Let $\Omega$ be an Eulerian regular CW-complex of dimension $d$. We say that $\Omega$ is \textit{\il{S}-shellable} (short for ``spherically shellable'') if either $\Omega=\{\emptyset\}$, or else we can linearly order the facets (open $d$-cells) of $\Omega$, say $\sigma_1,\hdots, \sigma_t$ such that for all $1\le i \le t$ the following conditions hold (where both \ $\bar{ }$ \ and \text{cl} denote closure).
\begin{enumerate}[a)]
    \item $\partial\overline{\sigma_1}$ is \il{S}-shellable (of dimension $d-1$)
    \item For $2\le i \le t-1$, let 
    $$\Gamma_{\sigma_i}=\text{cl}[\partial \overline{\sigma_i}-((\overline{\sigma_1}\cup \hdots\cup \overline{\sigma_{i-1}})\cap \overline{\sigma_i})].$$
\end{enumerate}
Thus $\Gamma_{\sigma_i}$ is the subcomplex of $\partial \overline{\sigma_i}$ generated by all $(d-1)$-cells of $\partial\overline{\sigma_i}$ which are not contained in the complex $\overline{\sigma_1}\cup \hdots\cup \overline{\sigma_{i-1}}$ generated by previous cells. Then we require that $\Gamma_{\sigma_i}$ is near-Eulerian of dimension $d-1$, and that the semisuspension  $\tilde\Sigma\Gamma_{\sigma_i}$ is \il{S}-shellable, with first facet of the shelling being the facet $\tau=\tau_i$ adjoined to $\Gamma_{\sigma_i}$ to obtain $\tilde\Sigma\Gamma_{\sigma_i}$.
\end{definition}

\begin{corollary}[Stanley~\cite{Stanley1994flag}]
The $\g{cd}$-index of every $S$-shellable sphere is nonnegative. 
\end{corollary}

\begin{proof}
Assume $\sigma_1, \sigma_2,\ldots,\sigma_t$ is an $S$-shelling of $\Omega$. If $\Omega=S^0$, then it has a \g{cd}-index of $\g{c}$. Thus, we may assume $\Omega$ has dimension at least $1$. Define $\Omega_i\coloneq \tilde{\Sigma}(\overline{\sigma}_1 \cup \cdots \cup \overline{\sigma}_i)$. Then $\Phi(\Omega_1)=\Phi(\partial \overline{\sigma}_1)\cdot \g{c}$ because $\Omega_1$ only has an initial coatom and a terminal coatom. Since $\partial\overline{\sigma}_1$ is assumed to be $S$-shellable, it has a nonnegative $\g{cd}$-index by induction. Then, when moving from $\Omega_{i-1}$ to $\Omega_i$, one adds all elements in $P_i$, the $i^{\text{th}}$ $S$-partition class, corresponding to the ordinary (or terminal) coatom $\sigma_i$. This is exactly the contribution of $\sigma_i$, witnessed by all paths of $T_\nu(\Omega)$ that include $P_i$, which is manifestly nonnegative. Hence the $\g{cd}$-index can only increase when moving from $\Omega_{i-1}$ to $\Omega_{i}$ and thus stays nonnegative. 
\end{proof}

\begin{figure}[h!]
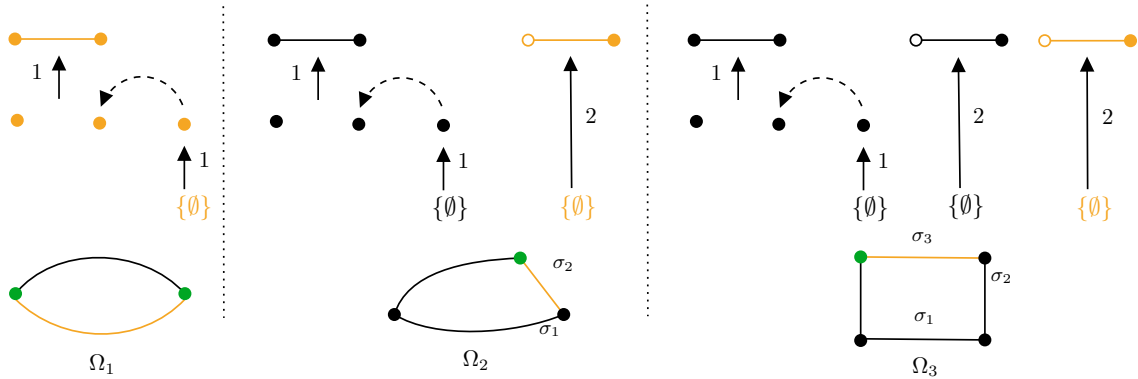

    \centering
\begin{center}
\adjustbox{scale = 0.85, center}{

\tikzset{every picture/.style={line width=0.75pt}} 



}
\end{center}
    \caption{The \g{cd}-index of each $\Omega_i$ for the square.}  \label{fig:Stanleysquare}
\end{figure}

In summary, Stanley's proof~\cite{Stanley1994flag} of the nonnegativity of the $\g{cd}$-index for all $S$-shellable spheres follows by adding a subtree to $T_\nu(\Omega_i)$ attached to the root $\Omega$ at each step. We illustrate this for the $\g{cd}$-index of the square in Figure~\ref{fig:Stanleysquare}.

\subsection{Shelling Components of the Simplex} 

There has been much study on the $\g{cd}$-index of the $d$-simplex and the $\g{cd}$-index of each of its shelling components. Mark Purtill~\cite{Purtill_1993} first established the nonnegativity of the $\g{cd}$-index of the $d$-simplex by studying Andr\'e permutations developed by Foata--Sch\"utzenberger~\cite{Foata}. Then, Stanley~\cite{Stanley1994flag} showed that the $\g{cd}$-index of any Eulerian simplicial poset depends only its $h$-vector and the $\g{cd}$-index of each shelling component of the simplex. Motivated by Purtill's result, Stanley conjectured that the $\g{cd}$-index of each shelling component of the simplex could be described via Andr\'e permutations in a similar manner. This conjecture was proven by G\'abor Hetyei~\cite{Hetyei1996}. Further extensions of relating Andr\'e permutations to the study of the $\g{cd}$-index can be found in Ehrenborg--Readdy~\cite{Ehrenborgcube}.

Alternatively, to study the \g{cd}-index of the $d$-simplex one can take the more algebraic approach due to Ehrenborg--Readdy~\cite{Ehrenborg1998}. The $\g{cd}$-index of the simplex can be computed by iterating the $\textbf{Pyr}$ operator, which acts as a derivation on $\mathbb{Z}\langle \g{c}, \g{d} \rangle$, and the \g{cd}-index of the shelling components of the simplex are computed via the derivation $\textbf{G}$. In this work, they establish similar algebraic formulas for the operations of prism, pyramid and bipyramid on polytopes. 

Given that there is much interest in the $\g{cd}$-index of the simplex and its shelling components, we take a moment to highlight a consequence of our combinatorial interpretation of the $\g{cd}$-index to the study of computing the $\g{cd}$-index of the $d$-simplex and its shelling components. Use $\hat{\Phi}_k^d$ to denote the $\g{cd}$-index of the shelling component with a restriction face of size $k$ of the boundary of the $d$-simplex. Also, define both $\hat{\Phi}_k^d[\g{c}]$ and $\hat{\Phi}_k^d[\g{d}]$ to be the parts of the $\g{cd}$-index of $\hat{\Phi}_k^d$ that are divisible on the right by $\g{c}$ and $\g{d}$ respectively. Naturally $\hat{\Phi}_k^d=\hat{\Phi}_k^d[\g{c}]+\hat{\Phi}_k^d[\g{d}]$. Furthermore, use $\Delta^{d}$ to denote the $d$-simplex and, in the statement of Corollary~\ref{coro:simplexrec}, use $\ast$ to denote the simplicial join of two simplicial complexes.

\begin{corollary}\label{coro:simplexrec}
    Let $d\geq 1$. Then for $0\leq k\leq d-1$, we obtain
    $$\hat{\Phi}_k^d[\g{c}]=\sum_{j=k}^{d-1}\hat{\Phi}_j^{d-1}\cdot \g{c}$$
and for $0< k <d$, we obtain 
$$\hat{\Phi}_k^d[\g{d}]=\sum_{j=0}^{d-2} h_j(\partial\Delta^{k-1} \ast \partial\Delta^{d-k-1})\hat{\Phi}_j^{d-2}\cdot \g{d}.$$
\end{corollary}

\begin{proof}
    Consider $\hat{\Phi}_k^d$. If $k=0$, then this shelling component is initial and is simply a $(d-1)$-simplex and we obtain $$\hat{\Phi}_0^d=\sum_{j=0}^{d-1}\hat{\Phi}_j^{d-1}\cdot \g{c}$$
    as desired. If $k=d$, then $\hat{\Phi}_{k}^d$ is terminal and has no $\g{cd}$-index contribution. If $0<k<d$, then we consider its $\g{cd}$-index as coming from an ordinary $S$-partition class. To compute the part of its $\g{cd}$-index that is right divisible by $\g{c}$, we consider $\tilde{\Sigma}\Gamma_k$, which is the ($d-1$)-sphere that is the result of merging the first $k$ facets in any shelling of the boundary of the $d$-simplex into one initial coatom $\tau$. The remaining contributions are given precisely by $\sum_{j=k}^{d-1}\hat{\Phi}_j^{d-1}\cdot \g{c}$. Meanwhile, $\partial\Gamma_k$ can be seen to be precisely $\partial\Delta^{k-1} \ast \partial\Delta^{d-k-1}$. Therefore, $\hat{\Phi}_k^d[\g{d}]$ is given by $\Phi(\partial\Delta^{k-1} \ast \partial\Delta^{d-k-1})\cdot \g{d}$. However, since $\partial\Delta^{k-1} \ast \partial\Delta^{d-k-1}$ is simplicial, we may apply Theorem~3.1 of Stanley~\cite{Stanley1994flag} to achieve the desired result.
\end{proof}

From this recursion, we can quickly recover a monotonicity result, a symmetry result, and a unimodality result concerning the $\g{cd}$-indices of shelling components of the simplex in Corollary~\ref{coro:simplexineqs}. For the unimodality result, we will need the following lemma, and the facts that the $h$-vector of the boundary of the $d$-simplex is the all one's vector of length $d+1$ and that $h$-polynomials are multiplicative over simplicial join.

\begin{lemma}\label{lem:polyprod}
    Let $d\geq 2$. For any fixed $0\leq \ell \leq d-2$, the sequence $\left\{\left(\sum_{i=0}^{k-1}t^{i}\right)\left(\sum_{j=0}^{d-k-1}t^{j}\right)[t^\ell]\right\}_{k=0,\ldots,d}$ is unimodal with respect to $k$, acquiring a maximum at $k=\lfloor d/2 \rfloor,\lceil d/2 \rceil$. 
\end{lemma}

\begin{proof}
    Because of the symmetry of the polynomial $\left(\sum_{i=0}^{k-1}t^{i}\right)\left(\sum_{j=0}^{d-k-1}t^{j}\right)$ when replacing $k$ with $d-k$, we will assume $k-1<d-k-1$, or simply $k<d-k$. Furthermore, our proof will follow by showing that if $k$ and $d-k$ differ by at least one, then replacing $k-1$ with $k$ and $d-k-1$ with $d-k-2$ will weakly increase every coefficient of the product. To this end, consider the following sequence of equalities:
    \begin{align*}
\left(\sum_{i=0}^{k}t^{i}\right)\left(\sum_{j=0}^{d-k-2}t^{j}\right)-\left(\sum_{i=0}^{k-1}t^{i}\right)\left(\sum_{j=0}^{d-k-1}t^{j}\right)&=t^{k}\left(\sum_{j=0}^{d-k-2}t^{j}\right)-t^{d-k-1}\left(\sum_{i=0}^{k-1}t^i\right)\\
&=\left(\sum_{j=k}^{d-2}t^{j}\right)-\left(\sum_{i=d-k-1}^{d-2}t^i\right).\\
\end{align*}
Since $k\leq d-k-1$, this difference is nonnegative as desired.
\end{proof}

\begin{corollary}\label{coro:simplexineqs}
    The following results hold for $d\geq 0$:
    \begin{enumerate}
        \item $\hat{\Phi}_0^d[\g{c}]\geq \hat{\Phi}_1^d[\g{c}] \geq \hat{\Phi}_2^d[\g{c}]\geq \cdots \geq \hat{\Phi}_d^d[\g{c}]$.
        \item $\hat{\Phi}_k^d[\g{d}]=\hat{\Phi}_{d-k}^d[\g{d}]$ for $0 \leq k \leq d$.
        \item $\hat{\Phi}_{\lfloor d/2 \rfloor }^d[\g{d}]\geq \cdots \geq \hat{\Phi}_2^d[\g{d}] \geq \hat{\Phi}_1^d[\g{d}] \geq \hat{\Phi}_0^d[\g{d}]=0$.
    \end{enumerate}  
\end{corollary}
\begin{proof} All three statements can be quickly verified for $d=0$ and $d=1$. Then, both (1) and (2) are immediate from Corollary~\ref{coro:simplexrec}. For part (3), recall that the $h$-polynomial of the boundary of the $k$-simplex is $1+t+t^2+\cdots+t^k$ and that $h$-polynomials are multiplicative over simplicial join. With this in mind, part (3) follows from Corollary~\ref{coro:simplexrec} and Lemma~\ref{lem:polyprod}.
    \end{proof}

\begin{figure}[h!]
\begin{subfigure}{.5\textwidth}
  \centering
\adjustbox{scale = 0.85, center}{
\begin{tikzcd}[sep=small]
	&&&& {\hat{\Phi}_4^4[\bf{c}]} & \cdots \\
	&&& {\hat{\Phi}_3^3[\bf{c}]} & {\hat{\Phi}_3^4[\bf{c}]} & \cdots \\
	&& {\hat{\Phi}_2^2[\bf{c}]} & {\hat{\Phi}_2^3[\bf{c}]} & {\hat{\Phi}_2^4[\bf{c}]} & \cdots \\
	& {\hat{\Phi}_1^1[\bf{c}]} & {\hat{\Phi}_1^2[\bf{c}]} & {\hat{\Phi}_1^3[\bf{c}]} & {\hat{\Phi}_1^4[\bf{c}]} & \cdots \\
	{\hat{\Phi}_0^0[\bf{c}]} & {\hat{\Phi}_0^1[\bf{c}]} & {\hat{\Phi}_0^2[\bf{c}]} & {\hat{\Phi}_0^3[\bf{c}]} & {\hat{\Phi}_0^4[\bf{c}]} & \cdots
	\arrow[hook, from=1-5, to=1-6]
	\arrow[hook, from=1-5, to=2-5]
	\arrow[hook, from=2-4, to=2-5]
	\arrow[hook, from=2-4, to=3-4]
	\arrow[hook, from=2-5, to=2-6]
	\arrow[hook, from=2-5, to=3-5]
	\arrow[hook, from=3-3, to=3-4]
	\arrow[hook, from=3-3, to=4-3]
	\arrow[hook, from=3-4, to=3-5]
	\arrow[hook, from=3-4, to=4-4]
	\arrow[hook, from=3-5, to=3-6]
	\arrow[hook, from=3-5, to=4-5]
	\arrow[hook, from=4-2, to=4-3]
	\arrow[hook, from=4-2, to=5-2]
	\arrow[hook, from=4-3, to=4-4]
	\arrow[hook, from=4-3, to=5-3]
	\arrow[hook, from=4-4, to=4-5]
	\arrow[hook, from=4-4, to=5-4]
	\arrow[hook, from=4-5, to=4-6]
	\arrow[hook, from=4-5, to=5-5]
	\arrow[hook, from=5-1, to=5-2]
	\arrow[hook, from=5-2, to=5-3]
	\arrow[hook, from=5-3, to=5-4]
	\arrow[hook, from=5-4, to=5-5]
	\arrow[hook, from=5-5, to=5-6]
\end{tikzcd}
}
  \caption{Inequalities for $\hat{\Phi}_k^d[\g{c}]$}
  \label{fig:sub1}
\end{subfigure}%
\begin{subfigure}{.5\textwidth}
  \centering
  \adjustbox{scale = 0.85, center}{
\begin{tikzcd}[sep=small]
	&&&& {\hat{\Phi}_4^4[\bf{d}]} & \cdots \\
	&&& {\hat{\Phi}_3^3[\bf{d}]} & {\hat{\Phi}_3^4[\bf{d}]} & \cdots \\
	&& {\hat{\Phi}_2^2[\bf{d}]} & {\hat{\Phi}_2^3[\bf{d}]} & {\hat{\Phi}_2^4[\bf{d}]} & \cdots \\
	& {\hat{\Phi}_1^1[\bf{d}]} & {\hat{\Phi}_1^2[\bf{d}]} & {\hat{\Phi}_1^3[\bf{d}]} & {\hat{\Phi}_1^4[\bf{d}]} & \cdots \\
	{\hat{\Phi}_0^0[\bf{d}]} & {\hat{\Phi}_0^1[\bf{d}]} & {\hat{\Phi}_0^2[\bf{d}]} & {\hat{\Phi}_0^3[\bf{d}]} & {\hat{\Phi}_0^4[\bf{d}]} & \cdots
	\arrow[hook, from=1-5, to=1-6]
	\arrow[hook, from=1-5, to=2-5]
	\arrow[hook, from=2-4, to=2-5]
	\arrow[hook, from=2-4, to=3-4]
	\arrow[hook, from=2-5, to=2-6]
	\arrow[hook, from=2-5, to=3-5]
	\arrow[hook, from=3-3, to=3-4]
	\arrow[hook, from=3-3, to=4-3]
	\arrow[hook, from=3-4, to=3-5]
	\arrow[equal, from=3-4, to=4-4]
	\arrow[hook, from=3-5, to=3-6]
	\arrow[hook, from=4-2, to=4-3]
	\arrow[hook, from=4-3, to=4-4]
	\arrow[hook, from=4-4, to=4-5]
	\arrow[hook', from=4-5, to=3-5]
	\arrow[hook, from=4-5, to=4-6]
	\arrow[equal, from=5-1, to=5-2]
	\arrow[equal, from=5-2, to=4-2]
	\arrow[equal, from=5-2, to=5-3]
	\arrow[hook', from=5-3, to=4-3]
	\arrow[equal, from=5-3, to=5-4]
	\arrow[hook', from=5-4, to=4-4]
	\arrow[equal, from=5-4, to=5-5]
	\arrow[hook', from=5-5, to=4-5]
	\arrow[hook, from=5-5, to=5-6]
\end{tikzcd}
}
  \caption{Inequalities for $\hat{\Phi}_k^d[\g{d}]$}
  \label{fig:sub2}
\end{subfigure}
\caption{Inequalities of Corollary~\ref{coro:simplexineqs}}
\label{fig:simplexeqs}
\end{figure}
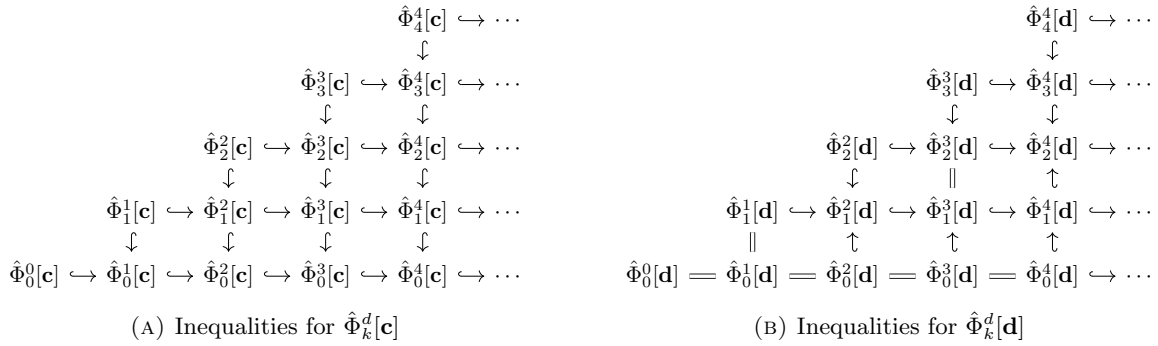

We shall use $\hookrightarrow$ to denote an inclusion on (multi)-sets of $\g{cd}$-words in Figure~\ref{fig:simplexeqs} where we provide illustrations for the inequalities of Corollary~\ref{coro:simplexineqs}.

Therefore, we have established another recurrence and multiple consequences for the shelling components of the simplex without the use of Andr\'e permutations. Next, we apply similar techniques to better understand the $\g{cd}$-index of the $d$-cube and its shelling components.

\subsection{Cubical Shelling Components}
In the case of the $d$-simplex, the shelling components are the upper intervals generated by a single element of the poset. Moreover, any linear order of the facets of a simplex is a shelling and an $S$-shelling. Hence, there is no confusion as to which type of shelling is needed. In the case of cubical spheres, it is also true that $C$-shellings are equivalent to $S$-shellings by Theorem~2.8 of Ehrenborg--Hetyei~\cite{EhrenborgHetyei2000}. Thus, there is again no need to distinguish between the two types of shellings.

Observe that neither a $d$-ball, $B$, nor a $d$-ball with its boundary removed, $B^\circ$, are Eulerian. Therefore, the $\g{cd}$-index need not exist a priori. However, we will take the definition of the $\g{cd}$-index of $B$ to be the $\g{cd}$-index of $\tilde{\Sigma}B$. As for $B^\circ$, we will consider the result of gluing a single $d$-cell $C$ to the boundary of $B$. By computing the $\g{cd}$-index of this resulting $d$-sphere, and subtracting $\Phi(\partial C)\cdot \g{c}$, we acquire $\Phi(B^\circ)$. 

Notice that in both cases, this definition of the $\g{cd}$-index can be understood as follows: Compute the flag $f$-polynomial of the poset in question. Transform this flag $f$-polynomial into the corresponding flag $h$-polynomial by the map $\g{a}\mapsto \g{a}-\g{b}$ and $\g{b}\mapsto \g{b}$. Replace the coefficient of every word of the form $\g{w}\g{b}$ with the coefficient of $\g{w}\g{a}$ in the flag $h$-polynomial of $\Omega$. In this way, a $\g{cd}$-index can be understood for any collection of $\g{ab}$-words of any flag $f$-polynomial, and when restricted to Eulerian posets or semi-Eulerian posets recovers the usual definitions of Bayer--Klapper~\cite{Bayer1991} and Juhnke--Samper--Venturello~\cite{semi}.

Alternatively, for the purpose of Proposition~\ref{prop:connsum}, one can use $\Phi(S)$ in place of $\Phi(B_1)+\Phi(B_2^\circ)$ where $S$ is the sphere obtained by gluing $B_1$ and $B_2$ along their combinatorially isomorphic boundaries.

\begin{proposition}\label{prop:connsum}
    Let $\Omega_1$ and $\Omega_2$ be CW $d$-spheres. Consider $\Omega_1 \# \Omega_2$ where the connected sum is formed by removing a $d$-ball $B_1$ of $\Omega_1$ and the interior of a $d$-ball $B_2$ where $\partial B_1$ is combinatorially isomorphic to $\partial B_2$. Then, the $\g{cd}$-index of the connected sum $\Omega_1 \# \Omega_2$ can be computed as follows: 
    $$\Phi(\Omega_1)+\Phi(\Omega_2)=\Phi(\Omega_1 \# \Omega_2)+\Phi(B_1)+\Phi(B_2^\circ).$$
\end{proposition}

\begin{proof}
    The proof follows directly from counting chains. Any element of $\Delta(\Omega_1) \setminus \Delta(\Omega_2)$ or $\Delta(\Omega_2)\setminus \Delta(\Omega_1)$ appears exactly once on each side. Any element of $ \Delta(\partial B_1)$ appears exactly twice on each side, once in each of $\Delta(\Omega_1)$ and $\Delta(\Omega_2$) on the left hand side and once in each of $\Delta(\Omega_1 \# \Omega_2)$ and $\Delta(B_1)$ on the right hand side. Additionally, any element of $\Delta(B_2)\setminus\Delta(\partial B_2)$ appears once on each side.   
\end{proof}

Using Proposition~\ref{prop:connsum}, we will establish an analog of the inequalities of Corollary~\ref{coro:simplexineqs} for cubical shelling components. Recall that a nonempty proper collection $\mathcal{C}$ of facets of the boundary of the $d$-cube forms an initial segment of a shelling if and only if there exists some facet $F\in \mathcal{C}$ whose antipode $-F$ is not an element of $\mathcal{C}$. Use $r$ to denote the number of facets present in $\mathcal{C}$ without their antipode and $s$ to denote the number of pairs of antipodal facets present in $\mathcal{C}$; we will say that this collection is of type $(r,s)$. Notice that all collections of facets of the boundary of the $d$-cube of type $(r,s)$ are combinatorially isomorphic as any two are related by a signed permutation of coordinates on the standard $d$-cube $[-1,1]^d$. Let $\mathcal{C}^d_{r,s}$ denote the shelling component of the boundary of the $(d+1)$-cube of type $(r,s)$. That is, for $r>0$, the shelling component $\mathcal{C}^{d}_{r,s}$ is the solid $d$-cube with a $(d-1)$-ball $B_{r,d-s-r}$ removed from its boundary, where $B_{r,d-s-r}$ is a ball generated by $r$ facets without their antipodal pairs and $d-s-r$ pairs of antipodal facets. Recall that $\mathcal{C}_{r,s}^d$ is a valid shelling component if and only if $r+s \leq d$ and either a) $r>0$ or b) $r=s=0$ or c) $r=0$ and $s=d$. Notice that the case $r=s=0$ corresponds to the terminal facet of the boundary of the $(d+1)$-cube while the case of $r=0$ and $s=d$ corresponds to the initial facet of the boundary of the $(d+1)$-cube.

\begin{corollary}\label{coro:cubeineqs}
The following inequalities hold whenever the relevant $\mathcal{C}_{r,s}^d$ form a shelling component: 
\begin{enumerate}
    \item $\Phi(\mathcal{C}^d_{r,s})[\g{c}]\leq \Phi(\mathcal{C}^d_{r+1,s})[\g{c}]$.
    \item $\Phi(\mathcal{C}^d_{r,s})[\g{c}]\leq \Phi(\mathcal{C}^d_{r-1,s+1})[\g{c}]$.
    \item $\Phi(\mathcal{C}^d_{r,s})[\g{d}]=\Phi(\mathcal{C}^d_{r,d-s-r})[\g{d}]$.
    \item $\Phi(\mathcal{C}^d_{r,s})[\g{d}]\leq \Phi(\mathcal{C}^d_{r+1,s})[\g{d}]$ and $r+2s\leq d-1$, with equality if and only if $r+2s=d-1$.
    \item $\Phi(\mathcal{C}^d_{r,s})[\g{d}]\leq \Phi(\mathcal{C}^d_{r-1,s+1})[\g{d}]$ and $r-1+2s\leq d-1$, with equality if and only if $r-1+2s=d-1$.
\end{enumerate}
    
\end{corollary}

\begin{proof}
Consider the cubical shelling component $\mathcal{C}_{r,s}^d$. Its \g{cd}-index can be computed recursively by Theorem~\ref{thm:combint}. In particular, $\Phi(\mathcal{C}_{r,s}^d)[\g{c}]$ is the \g{cd}-index of a terminal segment of a shelling of the boundary of the $d$-cube of type $(r,s)$, multiplied on the right by \g{c}. Recall that all shellings of the boundary of the $d$-cube are reversible. Hence, a set of facets forms an initial segment if and only if it forms a terminal segment.

(1) Observe that, by the existence of $\mathcal{C}^d_{r+1,s}$ as a shelling component, one can obtain a shelling of the boundary of the $d$-cube that finishes with a collection of facets such that there are $r+1$ facets without their antipodes and $s$ pairs of antipodal facets. Moreover, this terminal segment $T$ can be chosen so that the first facet $F$ of this segment satisfies $-F\not \in T$, by the existence of $\mathcal{C}^d_{r,s}$ as a shelling component. The difference $\Phi(\mathcal{C}_{r,s}^d)[\g{c}]-\Phi(\mathcal{C}^d_{r+1,s})[\g{c}]$ is precisely the contribution of $F$ in this shelling of the boundary of the $d$-cube multiplied on the right by \g{c}; thus it is nonnegative as desired.

(2) Observe that, by the existence of $\mathcal{C}^d_{r-1,s+1}$ as a shelling component, one can obtain a shelling of the boundary of the $d$-cube that finishes with a collection of facets such that there are $r-1$ facets without their antipodes and $s+1$ pairs of antipodal facets. Moreover, by the existence of $\mathcal{C}^d_{r,s}$ as a shelling component,  this terminal segment $T'$ can be chosen so that the first facet $F'$ of this segment satisfies $-F'\in T'$. The difference $\Phi(\mathcal{C}_{r-1,s+1}^d)[\g{c}]-\Phi(\mathcal{C}_{r,s}^d)[\g{c}]$ is precisely the contribution of $F'$ in this shelling of the boundary of the $d$-cube multiplied on the right by \g{c}; thus it is nonnegative as desired. 

(3) Use $\Gamma_{r,s}^d$ to denote the down closure of the elements of $\mathcal{C}_{r,s}^d$ in the face poset of the $d$-cube after deleting the single face of dimension $d$ from $\mathcal{C}_{r,s}^d$. Then, the statement follows immediately since $\partial \Gamma_{r,s}^d = \partial \Gamma_{r,d-s-r}^d$. 

(4) We will utilize Proposition~\ref{prop:connsum} because, when adding a facet to an initial segment of the boundary of a $d$-cube, the boundary of the resulting segment is the connected sum taken over the intersection of the newly added facet with all previous facets. Let $I$ denote an initial segment of a shelling of the boundary of the $d$-cube consisting of $r$ facets without their antipodal pairs and $s$ pairs of antipodal facets. By appending a facet $G$ whose antipode is missing, we obtain $\Gamma_{r+1,s}^d$. We apply Proposition~\ref{prop:connsum} with $\Omega_1$ equal to $\partial I$ and $\Omega_2=\partial G$. Furthermore, we take $B_1=B_2=G \cap I$. Notice that $G$ is a $(d-1)$-cube and has $2(d-1)$ facets. Thus, by Proposition~\ref{prop:connsum}, the desired change $\Phi(\mathcal{C}_{r+1,s}^d)[\g{d}]-\Phi(\mathcal{C}_{r,s}^d)[\g{d}]=\Phi(\partial\Gamma_{r+1,s}^d)\cdot \g{d}-\Phi(\partial\Gamma_{r,s}^d)\cdot \g{d}$ is given by $[\Phi(\partial G \# \partial I)-\Phi(I)]\cdot \g{d}=[\Phi(G)-\Phi(G\cap I)-\Phi((G \cap I)^\circ)]\cdot \g{d}$. Record the number of facets of $\partial G$ present in $G\cap I$ in a $0,1,2$ vector $\vec{v}$ of length $d-1$; each coordinate indicates the number of facets present from each corresponding antipodal pair. The condition $r+2s\leq d-1$ is equivalent to the condition that there are at least as many $0$ coordinates as there are $2$ coordinates in the vector $\vec{v}$. Let $\vec{w}$ be a vector obtained from $\vec{v}$ by fixing all $1$ coordinates and swapping the position of each coordinate equal to $2$ with the position of a coordinate equal to $0$. Next, we will shell the boundary of the $(d-1)$-cube so that two sets of facets each isomorphic to $G \cap I$ form an initial segment and a terminal segment of this shelling. That way, the difference $[\Phi(G)-\Phi(G\cap I)-\Phi((G \cap I)^\circ)]\cdot \g{d}$ will be witnessed by the \g{cd}-index contributions of the remaining facets of $G$ multiplied on the right by \g{d}. Form this desired shelling order by shelling the associated facets of $G$ corresponding to $\vec{v}$ and $\vec{w}$ as follows:
a) the coordinates equal to $1$ in $\vec{v}$, b) the coordinates equal to $2$ in $\vec{v}$, c) the coordinates that equal $0$ in both $\vec{v}$ and $\vec{w}$, d) the coordinates equal to $2$ in $\vec{w}$, e) the coordinates equal to $1$ in $\vec{w}$. Hence the difference $\Phi(\mathcal{C}^d_{r+1,s})[\g{d}]-\Phi(\mathcal{C}^d_{r,s})[\g{d}]$ is precisely the contributions of the facets added in step (c) as desired. To complete the proof, notice that set of facets added in step (c) is nonempty if and only if $r+2s<d-1$ and that any facet in step (c) is ordinary and thus has a nonzero \g{cd}-index.

(5) This argument is analogous to that of (4); one must simply re-count the number of facets that intersect the newly added $(d-1)$-cube $G'$. 
\end{proof}

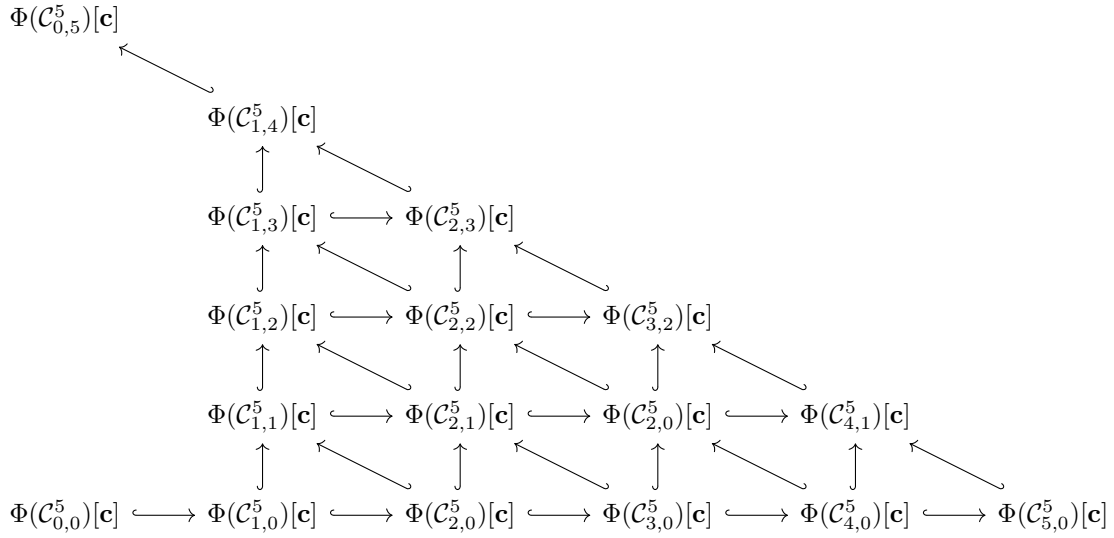
\begin{figure}[h!]
    \centering
    \begin{tikzcd}
	{\Phi(\mathcal{C}^5_{0,5})[\bf{c}]} &&&&& \\
	& {\Phi(\mathcal{C}_{1,4}^5)[\bf{c}]} \\
	& {\Phi(\mathcal{C}_{1,3}^5)[\bf{c}]} & {\Phi(\mathcal{C}_{2,3}^5)[\bf{c}]} \\
	& {\Phi(\mathcal{C}_{1,2}^5)[\bf{c}]} & {\Phi(\mathcal{C}_{2,2}^5)[\bf{c}]} & {\Phi(\mathcal{C}_{3,2}^5)[\bf{c}]} \\
	& {\Phi(\mathcal{C}_{1,1}^5)[\bf{c}]} & {\Phi(\mathcal{C}_{2,1}^5)[\bf{c}]} & {\Phi(\mathcal{C}_{2,0}^5)[\bf{c}]} & {\Phi(\mathcal{C}_{4,1}^5)[\bf{c}]} \\
	{\Phi(\mathcal{C}_{0,0}^5)[\bf{c}]} & {\Phi(\mathcal{C}_{1,0}^5)[\bf{c}]} & {\Phi(\mathcal{C}_{2,0}^5)[\bf{c}]} & {\Phi(\mathcal{C}_{3,0}^5)[\bf{c}]} & {\Phi(\mathcal{C}_{4,0}^5)[\bf{c}]} & {\Phi(\mathcal{C}_{5,0}^5)[\bf{c}]}
	\arrow[hook', from=2-2, to=1-1]
	\arrow[hook, from=3-2, to=2-2]
	\arrow[hook, from=3-2, to=3-3]
	\arrow[hook', from=3-3, to=2-2]
	\arrow[hook, from=4-2, to=3-2]
	\arrow[hook, from=4-2, to=4-3]
	\arrow[hook', from=4-3, to=3-2]
	\arrow[hook, from=4-3, to=3-3]
	\arrow[hook, from=4-3, to=4-4]
	\arrow[hook', from=4-4, to=3-3]
	\arrow[hook, from=5-2, to=4-2]
	\arrow[hook, from=5-2, to=5-3]
	\arrow[hook', from=5-3, to=4-2]
	\arrow[hook, from=5-3, to=4-3]
	\arrow[hook, from=5-3, to=5-4]
	\arrow[hook', from=5-4, to=4-3]
	\arrow[hook, from=5-4, to=4-4]
	\arrow[hook, from=5-4, to=5-5]
	\arrow[hook', from=5-5, to=4-4]
	\arrow[hook, from=6-1, to=6-2]
	\arrow[hook, from=6-2, to=5-2]
	\arrow[hook, from=6-2, to=6-3]
	\arrow[hook', from=6-3, to=5-2]
	\arrow[hook, from=6-3, to=5-3]
	\arrow[hook, from=6-3, to=6-4]
	\arrow[hook', from=6-4, to=5-3]
	\arrow[hook, from=6-4, to=5-4]
	\arrow[hook, from=6-4, to=6-5]
	\arrow[hook', from=6-5, to=5-4]
	\arrow[hook, from=6-5, to=5-5]
	\arrow[hook, from=6-5, to=6-6]
	\arrow[hook', from=6-6, to=5-5]
\end{tikzcd}
    \caption{Inequalities for $\Phi(\mathcal{C}^5_{r,s})[\g{c}]$}
    \label{fig:cubeineqc}
\end{figure}

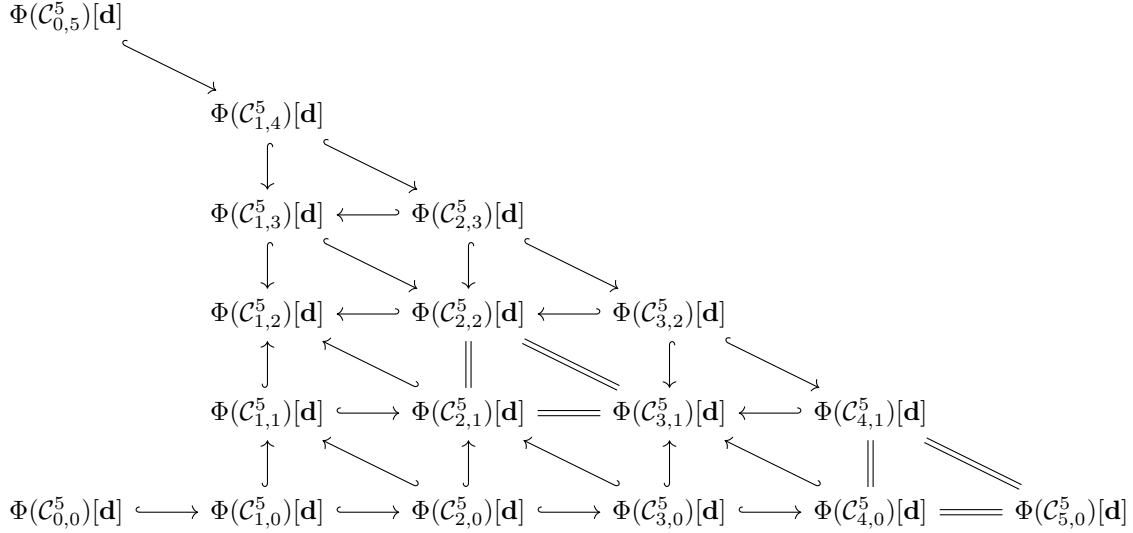
\begin{figure}[h!]
    \centering
\begin{tikzcd}
{\Phi(\mathcal{C}^5_{0,5})[\bf{d}]} &&&&& \\
	& {\Phi(\mathcal{C}_{1,4}^5)[\bf{d}]} \\
	& {\Phi(\mathcal{C}_{1,3}^5)[\bf{d}]} & {\Phi(\mathcal{C}_{2,3}^5)[\bf{d}]} \\
	& {\Phi(\mathcal{C}_{1,2}^5)[\bf{d}]} & {\Phi(\mathcal{C}_{2,2}^5)[\bf{d}]} & {\Phi(\mathcal{C}_{3,2}^5)[\bf{d}]} \\
	& {\Phi(\mathcal{C}_{1,1}^5)[\bf{d}]} & {\Phi(\mathcal{C}_{2,1}^5)[\bf{d}]} & {\Phi(\mathcal{C}_{3,1}^5)[\bf{d}]} & {\Phi(\mathcal{C}_{4,1}^5)[\bf{d}]} \\
	{\Phi(\mathcal{C}_{0,0}^5)[\bf{d}]} & {\Phi(\mathcal{C}_{1,0}^5)[\bf{d}]} & {\Phi(\mathcal{C}_{2,0}^5)[\bf{d}]} & {\Phi(\mathcal{C}_{3,0}^5)[\bf{d}]} & {\Phi(\mathcal{C}_{4,0}^5)[\bf{d}]} & {\Phi(\mathcal{C}_{5,0}^5)[\bf{d}]}
	\arrow[hook, from=1-1, to=2-2]
	\arrow[hook, from=2-2, to=3-2]
	\arrow[hook, from=2-2, to=3-3]
	\arrow[hook, from=3-2, to=4-2]
	\arrow[hook, from=3-2, to=4-3]
	\arrow[hook', from=3-3, to=3-2]
	\arrow[hook, from=3-3, to=4-3]
	\arrow[hook, from=3-3, to=4-4]
	\arrow[hook', from=4-3, to=4-2]
	\arrow[hook', from=4-4, to=4-3]
	\arrow[hook, from=4-4, to=5-4]
	\arrow[hook, from=4-4, to=5-5]
	\arrow[hook, from=5-2, to=4-2]
	\arrow[hook, from=5-2, to=5-3]
	\arrow[hook', from=5-3, to=4-2]
	\arrow[equals, nfold, from=5-3, to=4-3]
	\arrow[equals, nfold, from=5-3, to=5-4]
	\arrow[equals, nfold, from=5-4, to=4-3]
	\arrow[hook', from=5-5, to=5-4]
	\arrow[equals, nfold, from=5-5, to=6-6]
	\arrow[hook, from=6-1, to=6-2]
	\arrow[hook, from=6-2, to=5-2]
	\arrow[hook, from=6-2, to=6-3]
	\arrow[hook', from=6-3, to=5-2]
	\arrow[hook, from=6-3, to=5-3]
	\arrow[hook, from=6-3, to=6-4]
	\arrow[hook', from=6-4, to=5-3]
	\arrow[hook, from=6-4, to=5-4]
	\arrow[hook, from=6-4, to=6-5]
	\arrow[hook', from=6-5, to=5-4]
	\arrow[equals, nfold, from=6-5, to=5-5]
	\arrow[equals, nfold, from=6-5, to=6-6]
\end{tikzcd}
    \caption{Inequalities for $\Phi(\mathcal{C}^5_{r,s})[\g{d}]$}
    \label{fig:cubeineqsd}
\end{figure}

Again, using $\hookrightarrow$ to denote an inclusion on (multi)-sets of $\g{cd}$-words, we obtain illustrations for the inequalities of Corollary~\ref{coro:cubeineqs} in Figure~\ref{fig:cubeineqc} and Figure~\ref{fig:cubeineqsd}

\subsection{Strict Positivity of $\g{cd}$-indices}

The fact that all coefficients of the $\g{cd}$-index of any positive dimensional polytope are strictly positive has been known among experts. It is a consequence of the remark ``the coefficients of the $\g{cd}$-index of the simplex are all positive" made in the introduction of Billera--Ehrenborg \cite{Billera2000}, which can verified by iterating the Pyramid operator of Ehrenborg--Readdy \cite{Ehrenborg1998}. Then, when combined with the main result of Billera--Ehrenborg \cite{Billera2000} that the $\g{cd}$-index of any $d$-polytope is minimized on the $d$-simplex, one achieves strict positivity for all polytopes. We provide an alternative proof that follows immediately from Theorem~\ref{thm:combint}. 

\begin{corollary}\label{coro:strictpospoly}
    For all $d\geq 1$, the $\g{cd}$-index of every $d$-polytope is strictly positive on every $\g{cd}$-word. 
\end{corollary}

\begin{proof}
We shall proceed via induction on dimension. The $\g{cd}$-index of a line segment is $\g{c}$ and the $\g{cd}$-index of an $n$-gon is $\g{c}^2+(n-2)\g{d}$. For $d\geq 3$, we know that any $d$-polytope has at least $3$ facets. Take a line shelling and consider the initial facet $\sigma_1$ and the ordinary facet $\sigma_2$. By induction, the summands $\Phi(\partial \sigma_1)\cdot \g{c}$ and $\Phi(\partial \Gamma_2)\cdot \g{d}$ contain every $\g{cd}$-word, since $\partial \Gamma_2$ is polytopal as noted after Definition~2.1 of Stanley~\cite{Stanley1994flag}. 
\end{proof}

While the previous theorem holds for all polytopes, the same is not true for its shelling components. This can be seen even for the shelling components of the simplex by looking at the computations provided at the end of Stanley~\cite{Stanley1994flag}. This leads us to ask the following question: 

\begin{question}
    Which $\g{cd}$-monomials are supported in $\hat{\Phi}_i^d$ for each $0\leq i \leq d$?
\end{question}
The answer for $i=0$ follows immediately from Corollary~\ref{coro:strictpospoly}: every $\g{cd}$-monomial ending in $\g{c}$. Furthermore, for $i=d$, the shelling component is terminal and so no $\g{cd}$-monomials are supported. Using Corollary~\ref{coro:simplexrec}, the following results can be quickly derived by induction. 

\begin{corollary}
    $\hat{\Phi}_1^d$ is supported on all $\g{cd}$-monomials except $\g{c}^d$ for all $d\geq 2$.
\end{corollary}

\begin{corollary}
    $\hat{\Phi}_2^d$ is supported on all $\g{cd}$-monomials except $\g{c}^d$ and $\g{d}\g{c}^{d-2}$ for all $d\geq 3$.
\end{corollary}

Furthermore, because $\partial \Gamma_\sigma$ of an ordinary shelling component of the simplex is polytopal, we find that every \g{cd}-word ending in $\g{d}$ is supported in every $\hat{\Phi}_i^d$ for $0<i<d$ by Corollary~\ref{coro:strictpospoly}. Using the supports of the two shelling components $\hat{\Phi}_1^d$ and $\hat{\Phi}_2^d$, we will derive the proofs of Corollary~\ref{coro:unbdonone} and Corollary~\ref{coro:unbdontwo}.

Bayer~\cite{Bayer2001} established the suprema and infima of the coefficients of each $\g{cd}$-monomial as a function of all Eulerian posets as follows in Proposition~\ref{thm:bayerconstruct}.

\begin{proposition}\label{thm:bayerconstruct}(Bayer~\cite{Bayer2001})
For the following $\g{cd}$-words \g{w} as a function of Eulerian posets:
\begin{enumerate}
    \item The coefficients of the following have a greatest lower bound $0$ and have no upper bound:
    \begin{enumerate}
        \item $\g{c}^i\g{d}\g{c}^j$ with $\min\{i,j\} \leq 1$;
        \item $\g{c}^i\g{d}\g{c}\g{d}\cdots\g{c}\g{d}\g{c}^j$, with at least two $\g{d}$'s alternating with $\g{c}$'s, where $i$ and $j$ are unrestricted.
    \end{enumerate}
    \item The coefficient of $\g{c}^d$ of every Eulerian poset is $1$.
    \item For all other $\g{cd}$-words \g{w}, the coefficient of \g{w} as a function of Eulerian posets has neither a lower nor an upper bound.
\end{enumerate}
\end{proposition}

From Karu~\cite{Karu_2006}, we know that the $\g{cd}$-index of every Gorenstein* poset is nonnegative. Because the regular CW $d$-sphere with $2$ cells in every dimension has a $\g{cd}$-index of $\g{c}^{d+1}$, this lower bound is tight for all words not equal to $\g{c}^{d+1}$. In what follows, we will give a simple construction of a family of spheres that will show that the coefficient of every other $\g{cd}$-monomial can be arbitrarily large. Then, we will present a similar construction of a family of $d$-spheres on $d+2$ vertices that will have an unbounded $\g{cd}$-index for all $\g{cd}$-words except for the two monomials $\g{c}^{d+1}$ and $\g{d}\g{c}^{d-1}$. Since the coefficient of $\g{d}\g{c}^{d-1}$ equals the number of vertices minus two, this coefficient cannot be taken to be arbitrarily large on a fixed number of vertices. However, on a fixed number of vertices, every other monomial can be taken to be arbitrarily large. 

\begin{corollary}\label{coro:unbdonone}
    For every $\g{cd}$-polynomial $P(\g{c},\g{d})=\sum_{\g{w}\neq c^{d+1}}c_\g{w} \g{w}$, there is a subdivision of the boundary of the $(d+1)$-simplex whose $\g{cd}$-index is coefficient-wise greater than $P(\g{c},\g{d})$ for every $d\geq 2$.
\end{corollary}

\begin{proof}
    Start with the boundary of the $(d+1)$-simplex. On the second facet of the shelling, $F_2$ with restriction face $v$, place a vertex $v'$ and cone over the facet of $F_2$ not containing $v$. This operation subdivides the facet $F_2$, while adding an additional $S$-partition class of the form $\hat{\Phi}_1^d$. Indeed, $F_2$ is a $d$-ball, and another copy of $F_2$ is glued inside along a $(d-1)$-ball on its boundary. Since $\hat{\Phi}_1^d$ is supported on all monomials except $\g{c}^{d+1}$ by Corollary~\ref{coro:unbdonone}, we may repeat this operation more than $\text{max}_{\g{w}} \g{c}_{\g{w}}$ many times in order to exceed the coefficients of $P(\g{c},\g{d})$.
\end{proof}

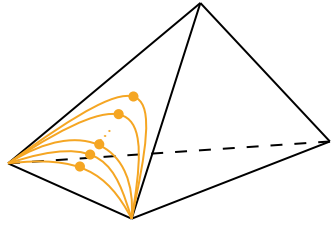
\begin{figure}[h!]
    \centering
    \begin{center}
\tikzset{every picture/.style={line width=0.75pt}} 

\begin{tikzpicture}[x=0.75pt,y=0.75pt,yscale=-1,xscale=1]

\draw    (301,399.84) -- (400.5,361.09) ;
\draw    (239,372.09) -- (301,399.84) ;
\draw  [dash pattern={on 4.5pt off 4.5pt}]  (239,372.09) -- (400.5,361.09) ;
\draw    (335.5,291.59) -- (301,399.84) ;
\draw    (335.5,291.59) -- (400.5,361.09) ;
\draw    (335.5,291.59) -- (239,372.09) ;
\draw [color={rgb, 255:red, 245; green, 166; blue, 35 }  ,draw opacity=1 ]   (239,372.09) .. controls (287.25,356.09) and (292.75,373.59) .. (301,399.84) ;
\draw [color={rgb, 255:red, 245; green, 166; blue, 35 }  ,draw opacity=1 ]   (239,372.09) .. controls (277.75,347.59) and (303.75,364.09) .. (301,399.84) ;
\draw [color={rgb, 255:red, 245; green, 166; blue, 35 }  ,draw opacity=1 ]   (239,372.09) .. controls (306.5,327.59) and (310,346.59) .. (301,399.84) ;
\draw [color={rgb, 255:red, 245; green, 166; blue, 35 }  ,draw opacity=1 ]   (239,372.09) .. controls (313.5,318.09) and (316,329.59) .. (301,399.84) ;
\draw [color={rgb, 255:red, 245; green, 166; blue, 35 }  ,draw opacity=1 ]   (239,372.09) .. controls (274.5,368.09) and (288.5,376.09) .. (301,399.84) ;
\draw  [color={rgb, 255:red, 245; green, 166; blue, 35 }  ,draw opacity=1 ][fill={rgb, 255:red, 245; green, 166; blue, 35 }  ,fill opacity=1 ] (299.86,338.52) .. controls (299.86,337.39) and (300.77,336.48) .. (301.89,336.48) .. controls (303.02,336.48) and (303.93,337.39) .. (303.93,338.52) .. controls (303.93,339.64) and (303.02,340.55) .. (301.89,340.55) .. controls (300.77,340.55) and (299.86,339.64) .. (299.86,338.52) -- cycle ;
\draw  [color={rgb, 255:red, 245; green, 166; blue, 35 }  ,draw opacity=1 ][fill={rgb, 255:red, 245; green, 166; blue, 35 }  ,fill opacity=1 ] (292.43,347.37) .. controls (292.43,346.25) and (293.34,345.34) .. (294.46,345.34) .. controls (295.59,345.34) and (296.5,346.25) .. (296.5,347.37) .. controls (296.5,348.5) and (295.59,349.41) .. (294.46,349.41) .. controls (293.34,349.41) and (292.43,348.5) .. (292.43,347.37) -- cycle ;
\draw  [color={rgb, 255:red, 245; green, 166; blue, 35 }  ,draw opacity=1 ][fill={rgb, 255:red, 245; green, 166; blue, 35 }  ,fill opacity=1 ] (282.71,362.52) .. controls (282.71,361.39) and (283.63,360.48) .. (284.75,360.48) .. controls (285.87,360.48) and (286.78,361.39) .. (286.78,362.52) .. controls (286.78,363.64) and (285.87,364.55) .. (284.75,364.55) .. controls (283.63,364.55) and (282.71,363.64) .. (282.71,362.52) -- cycle ;
\draw  [color={rgb, 255:red, 245; green, 166; blue, 35 }  ,draw opacity=1 ][fill={rgb, 255:red, 245; green, 166; blue, 35 }  ,fill opacity=1 ] (278.14,367.66) .. controls (278.14,366.53) and (279.05,365.62) .. (280.18,365.62) .. controls (281.3,365.62) and (282.21,366.53) .. (282.21,367.66) .. controls (282.21,368.78) and (281.3,369.69) .. (280.18,369.69) .. controls (279.05,369.69) and (278.14,368.78) .. (278.14,367.66) -- cycle ;
\draw  [color={rgb, 255:red, 245; green, 166; blue, 35 }  ,draw opacity=1 ][fill={rgb, 255:red, 245; green, 166; blue, 35 }  ,fill opacity=1 ] (273,373.66) .. controls (273,372.53) and (273.91,371.62) .. (275.04,371.62) .. controls (276.16,371.62) and (277.07,372.53) .. (277.07,373.66) .. controls (277.07,374.78) and (276.16,375.69) .. (275.04,375.69) .. controls (273.91,375.69) and (273,374.78) .. (273,373.66) -- cycle ;

\draw (281.92,361.34) node [anchor=north west][inner sep=0.75pt]  [font=\scriptsize,color={rgb, 255:red, 245; green, 166; blue, 35 }  ,opacity=1 ,rotate=-313.95]  {$...$};
\end{tikzpicture}
\end{center}
    \caption{Construction of Corollary~\ref{coro:unbdonone}}
\label{fig:unbdonone}
\end{figure}

\begin{corollary}\label{coro:unbdontwo}
    For every $\g{cd}$-polynomial $P(\g{c},\g{d})=\sum_{\g{w}\neq \g{c}^{d+1},\g{d}\g{c}^{d-1}}c_\g{w} \g{w}$, there is a subdivision of the boundary of the $(d+1)$-simplex whose $\g{cd}$-index is coefficient-wise greater than $P(\g{c},\g{d})$ for every $d\geq 3$ without adding any additional vertices.
\end{corollary}

\begin{proof}
      Start with the boundary of the $(d+1)$-simplex. Consider the shelling component $[e,\Delta_{d+1})$, where $e$ is an edge. The $\g{cd}$-index of this shelling component is supported on all monomials except $\g{c}^{d+1}$ and $\g{d}\g{c}^{d-1}$ by Corollary~\ref{coro:unbdontwo} and is topologically a $d$-ball with a $(d-1)$-ball removed from its boundary. By gluing in sufficiently many copies of this shelling component within $F_3$ along its intersection with $F_1$ and $F_2$, one creates a subdivision of the boundary of the $(d+1)$-simplex with the desired properties.
\end{proof}

\begin{figure}[h!]
    \centering
    \begin{center}
    
\tikzset{every picture/.style={line width=0.75pt}} 

\begin{tikzpicture}[x=0.75pt,y=0.75pt,yscale=-1,xscale=1]

\draw    (301,399.84) -- (400.5,361.09) ;
\draw    (239,372.09) -- (301,399.84) ;
\draw  [dash pattern={on 4.5pt off 4.5pt}]  (239,372.09) -- (400.5,361.09) ;
\draw    (335.5,291.59) -- (301,399.84) ;
\draw    (335.5,291.59) -- (400.5,361.09) ;
\draw    (335.5,291.59) -- (239,372.09) ;
\draw [color={rgb, 255:red, 245; green, 166; blue, 35 }  ,draw opacity=1 ]   (239,372.09) .. controls (287.25,356.09) and (292.75,373.59) .. (301,399.84) ;
\draw [color={rgb, 255:red, 245; green, 166; blue, 35 }  ,draw opacity=1 ]   (239,372.09) .. controls (277.75,347.59) and (303.75,364.09) .. (301,399.84) ;
\draw [color={rgb, 255:red, 245; green, 166; blue, 35 }  ,draw opacity=1 ]   (239,372.09) .. controls (306.5,327.59) and (310,346.59) .. (301,399.84) ;
\draw [color={rgb, 255:red, 245; green, 166; blue, 35 }  ,draw opacity=1 ]   (239,372.09) .. controls (313.5,318.09) and (316,329.59) .. (301,399.84) ;
\draw [color={rgb, 255:red, 245; green, 166; blue, 35 }  ,draw opacity=1 ]   (239,372.09) .. controls (274.5,368.09) and (288.5,376.09) .. (301,399.84) ;

\draw (281.92,361.34) node [anchor=north west][inner sep=0.75pt]  [font=\scriptsize,color={rgb, 255:red, 245; green, 166; blue, 35 }  ,opacity=1 ,rotate=-313.95]  {$...$};

\end{tikzpicture}

\end{center}
    \caption{Construction of Corollary~\ref{coro:unbdontwo}}
    \label{fig:unbdontwo}
\end{figure}
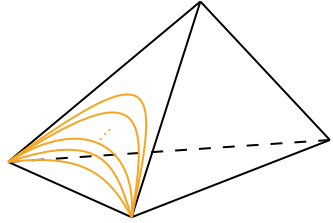

Next, we prove a new minor result regarding Eulerian posets with a free $\mathbb{Z}/2$-action. That is, a free involution $\tau$ on every element of the face poset that preserves the poset relation and the rank. That is, a) if $p \leq q$, then $\tau p \leq \tau q$ and b) if $\rho(p)=k$, then $\rho(\tau p)=k$. An example is the face poset of any centrally symmetric polytope with the antipodal action. 

\begin{proposition}\label{prop:cseven}
    Let $\Omega$ be an Eulerian poset of rank $d+1$ equipped with a free $\mathbb{Z}/2$-action. Then, the coefficient of every $\g{cd}$-word of $\Omega$ except for $\g{c}^{d+1}$ must be even.
\end{proposition}

\begin{proof}
    Every coefficient of the $\g{cd}$-index of $\Omega$ is a linear combination of the entries of the flag $f$-vector of $\Omega$. Observe that $f_K(\Omega)$ is even for every $K \neq \emptyset$ and $f_\emptyset(\Omega)=1$. Thus, the parity of each coefficient of the $\g{cd}$-index of $\Omega$ is determined by that of any single Eulerian poset with a free $\mathbb{Z}/2$-action. Since, the regular CW $d$-sphere with $2$ cells in every dimension has a $\g{cd}$-index of $\g{c}^{d+1}$, the result follows. 
\end{proof}

\begin{corollary}
    The $\g{cd}$-index of every centrally symmetric $(d+1)$-polytope has even coefficients for all $\g{cd}$-words not equal to $\g{c}^{d+1}$.
\end{corollary}

\subsection{Joins}
It was quickly observed in Lemma~1.1 of Stanley~\cite{Stanley1994flag} that the $\g{cd}$-index is multiplicative over poset join. Furthermore, it is remarked in Ehrenborg-Karu~\cite{Ehrenborg-Karu2007} that the join of two Gorenstein* posets is again Gorenstein*. However, it is not obvious that $S$-partitionability or $S$-shellability are preserved under poset join. In this subsection, we prove that both $S$-partitionability and $S$-shellability are preserved under poset join. While not immediate, these statements will follow from a careful unraveling of their definitions. Therefore, we will present only the proof for $S$-shellable spheres as the proof for $S$-partitionable posets is analogous.

\begin{theorem}\label{thm:Sshelljoin}
The (poset) join of two $S$-shellable spheres is $S$-shellable. 
\end{theorem}

\begin{proof}
Let $\Omega_1$ and $\Omega_2$ be two $S$-shellable spheres. We will prove that $\Omega_1 \ast \Omega_2$ is $S$-shellable via induction on the rank of $\Omega_2$. If $\Omega_2$ is the poset with only the distinct elements $\hat{0}$ and $\hat{1}$, then the claim quickly follows because $\Omega_1 \ast \Omega_2=\Omega_1$. Therefore, assume $\Omega_2$ has larger rank, assume the inductive hypothesis and consider $\Omega_1 \ast \Omega_2$. We shall show that the same $S$-shelling order of the facets of $\Omega_2$ will induce an $S$-shelling order on $\Omega_1 \ast \Omega_2$.  Declare the following shelling order:
\begin{itemize}
\item $\sigma_{1}^{\Omega_1 \ast \Omega_2}\coloneq \sigma_1^{\Omega_2}$ where $\sigma_1^{\Omega_2}$ was the initial coatom of $\Omega_2$, now viewed as a coatom of $\Omega_1 \ast \Omega_2$. 
\item $\sigma_i^{\Omega_1 \ast \Omega_2} \coloneq \sigma_i^{\Omega_2}$ for all $2\leq i \leq t-1$. 
\item $\sigma_t^{\Omega_1 \ast \Omega_2} \coloneq\sigma_t^{\Omega_2}$.
\end{itemize}
Next, we check that this total order on facets is indeed an $S$-shelling order by analyzing $\Gamma_i$ for each $i$. Observe that $\Gamma_i^{\Omega_1 \ast \Omega_2}=\Omega_1 \ast \Gamma_{i}^{\Omega_2}$ and $\partial \Gamma_i^{\Omega_2 \ast \Omega_1}=\Omega_1 \ast \partial \Gamma_i^{\Omega_2}$ for all $1<i<t$ by the definition of poset join. Furthermore, for every $1<i<t$, we have $\tilde{\Sigma}\Gamma_i^{\Omega_1 \ast \Omega_2}=\Omega_1 \ast (\tilde{\Sigma}\Gamma_i^{\Omega_2})$. Next, we verify the recursive properties of Definition~\ref{defn:Sshell} on each $\Gamma_i$.
\begin{enumerate}
\item $\overline{\sigma}_{1}^{\Omega_1 \ast \Omega_2}=\Omega_1 \ast \overline{\sigma}_1^{\Omega_2}$ and thus $\partial\overline{\sigma_1}^{\Omega_1 \ast \Omega_2}=\Omega_1\ast \partial\overline{\sigma}_1^{\Omega_2}$. Since $\Omega_2$ is $S$-shellable, we have that $\partial\overline{\sigma}_1^{\Omega_2}$ is $S$-shellable of rank one lower. Therefore, by induction, we find that $\partial \overline{\sigma}_1^{\Omega_1 \ast \Omega_2}$ is $S$-shellable.  
\item Since $\tilde{\Sigma}\Gamma_i^{\Omega_1 \ast \Omega_2}=\Omega_1 \ast (\tilde{\Sigma}\Gamma_i^{\Omega_2})$, we have by induction that $\tilde{\Sigma}\Gamma_i^{\Omega_1 \ast \Omega_2}$ is $S$-shellable of rank one lower. Furthermore, because the initial coatom of $\Omega_1 \ast \Omega_2$ corresponds to the initial coatom of $\Omega_2$, we find that the newly added coatom $\tau_i^{\Omega_1 \ast \Omega_2}$ is the initial coatom in $\tilde{\Sigma}\Gamma_i^{\Omega_1 \ast \Omega_2}$.
\end{enumerate}
This establishes the necessary inductive properties of Definition~\ref{defn:Sshell} as desired.  
\end{proof}

The following corollary immediately follows from Theorem~\ref{thm:Sshelljoin} via induction.

\begin{corollary}
The join of finitely many $S$-shellable spheres is an $S$-shellable sphere. 
\end{corollary}

\begin{remark}
    An analogous proof also shows that the join of any two $S$-partitionable posets is $S$-partitionable, and thus the join of finitely many $S$-partitionable posets is $S$-partitionable.
\end{remark}

\subsection{Interpretations of formulas from Ehrenborg and Readdy} In their 1998 paper, Ehrenborg and Readdy \cite{Ehrenborg1998} provided explicit formulas describing how the $\g{cd}$-index of a polytope changes under the effect of the prism, pyramid and bipyramid operations. They describe these operations in the generality of posets and give the following theorem: 
\begin{proposition}[Ehrenborg--Readdy~\cite{Ehrenborg1998}]
Let $P$ be an Eulerian poset. Then 
\begin{itemize}
    \item $\Phi(\operatorname{Prism}(P))=\Phi(P)\cdot\g{c}+D(\Phi(P))$,
    \item $\Phi(\operatorname{Bipyr}(P))=\g{c}\cdot\Phi(P)+D(\Phi(P))$,
    \item $2\Phi(\operatorname{Pyr}(P))=\Phi(P)\cdot \g{c}+\g{c}\cdot\Phi(P)+D(\Phi(P))$.
\end{itemize}
\end{proposition}
\noindent
Here $D$ is a derivation on $\g{cd}$-words with $D(\g{c})=2\g{d}$ and $D(\g{d})=\g{cd}+\g{dc}$. That is, $D$ is a linear map, and is recursively defined on $\g{cd}$-words via $D(\g{w}\g{c})=2\g{w}\g{d}+D(\g{w})\g{c}$ and $D(\g{w}\g{d})=\g{w}(\g{cd}+\g{dc})+D(\g{w})\g{d}$. In this short subsection, we take a moment to highlight how these formulas can be visualized via Theorem~\ref{thm:combint}.  

\begin{figure}[h!]
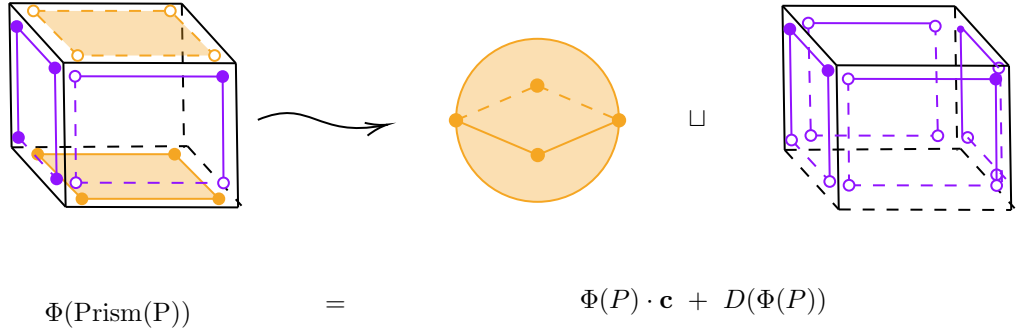

    \centering
    \begin{center}

\tikzset{every picture/.style={line width=0.75pt}} 



\end{center}
    \caption{Prism Operation}
    \label{fig:prism}
\end{figure}

We will begin with the prism operation, which takes a polytope $P$ to the polytope $P \times [0,1]$. Take a line shelling that begins with the facet $P \times \{0\}$ and ends with the facet $P \times \{1\}$. This can be done by Corollary~8.3 of Ziegler~\cite{Ziegler_1995}. The initial facet will contribute $\Phi(P)\cdot \g{c}$. The terminal facet will not contribute any $\g{cd}$-words; we include it with the initial facet to partition the faces into two Eulerian objects. Thus, the contributions of all facets of the form $F \times [0,1]$, where $F$ is a facet of $P$, will make up the $D(\Phi(P))$ term, while the initial and terminal facet together contribute $\Phi(P)\cdot \g{c}$. In Figure~\ref{fig:prism}, we visualize the contributions of the initial and terminal facet by gluing them together along their boundaries, which both equal $\partial P$, alongside the contributions of the facets of the form $F \times [0,1]$.

Observe that the prism operation takes a polytope $P$ to $P \times [0,1]$. Given two polytopes $P$ and $Q$, it is thus intriguing to ask for a formula for the \g{cd}-index of $P \times Q$ in terms of the \g{cd}-indices of $P$ and $Q$ individually. Phrased more generally, we ask the following:

\begin{question}
    Given two Eulerian posets $P$ and $Q$, can one describe the $\g{cd}$-index of $P \times Q$ in terms of the $\g{cd}$-index of $P$ and the $\g{cd}$-index of $Q$?
\end{question}

As for the Bipyramid operation, this will follow from the short dualization argument of Ehrenborg--Readdy~\cite{Ehrenborg1998} because $\text{Bipyr}(P)=\text{Prism}(P^*)^*$ and the effect of dualization on the \g{cd}-index is a reversal of all \g{cd}-monomials. 

\begin{figure}[h!]
    \centering
    \begin{center}
\tikzset{every picture/.style={line width=0.75pt}} 

\begin{tikzpicture}[x=0.75pt,y=0.75pt,yscale=-1,xscale=1]

\draw    (101,400.34) -- (142.13,421.59) ;
\draw  [dash pattern={on 4.5pt off 4.5pt}]  (152.88,389.34) -- (192.88,411.59) ;
\draw  [dash pattern={on 4.5pt off 4.5pt}]  (101,400.34) -- (152.88,389.34) ;
\draw    (142.13,421.59) -- (192.88,411.59) ;
\draw    (150.6,330.14) -- (192.88,411.59) ;
\draw    (150.6,330.14) -- (101,400.34) ;
\draw  [color={rgb, 255:red, 245; green, 166; blue, 35 }  ,draw opacity=1 ][fill={rgb, 255:red, 245; green, 166; blue, 35 }  ,fill opacity=0.5 ] (110.7,402.83) -- (143.21,418.85) -- (185.49,410.39) -- (152.98,394.36) -- cycle ;
\draw  [dash pattern={on 4.5pt off 4.5pt}]  (150.6,330.14) -- (152.88,389.34) ;
\draw  [color={rgb, 255:red, 245; green, 166; blue, 35 }  ,draw opacity=1 ][fill={rgb, 255:red, 245; green, 166; blue, 35 }  ,fill opacity=1 ] (109.85,402.34) .. controls (109.85,401.61) and (110.45,401.01) .. (111.18,401.01) .. controls (111.91,401.01) and (112.51,401.61) .. (112.51,402.34) .. controls (112.51,403.08) and (111.91,403.67) .. (111.18,403.67) .. controls (110.45,403.67) and (109.85,403.08) .. (109.85,402.34) -- cycle ;
\draw  [color={rgb, 255:red, 245; green, 166; blue, 35 }  ,draw opacity=1 ][fill={rgb, 255:red, 245; green, 166; blue, 35 }  ,fill opacity=1 ] (151.65,394.36) .. controls (151.65,393.63) and (152.24,393.03) .. (152.98,393.03) .. controls (153.71,393.03) and (154.3,393.63) .. (154.3,394.36) .. controls (154.3,395.09) and (153.71,395.69) .. (152.98,395.69) .. controls (152.24,395.69) and (151.65,395.09) .. (151.65,394.36) -- cycle ;
\draw  [color={rgb, 255:red, 245; green, 166; blue, 35 }  ,draw opacity=1 ][fill={rgb, 255:red, 245; green, 166; blue, 35 }  ,fill opacity=1 ] (184.16,410.39) .. controls (184.16,409.65) and (184.76,409.06) .. (185.49,409.06) .. controls (186.22,409.06) and (186.82,409.65) .. (186.82,410.39) .. controls (186.82,411.12) and (186.22,411.71) .. (185.49,411.71) .. controls (184.76,411.71) and (184.16,411.12) .. (184.16,410.39) -- cycle ;
\draw  [color={rgb, 255:red, 245; green, 166; blue, 35 }  ,draw opacity=1 ][fill={rgb, 255:red, 245; green, 166; blue, 35 }  ,fill opacity=1 ] (141.33,418.22) .. controls (141.33,417.49) and (141.92,416.89) .. (142.66,416.89) .. controls (143.39,416.89) and (143.98,417.49) .. (143.98,418.22) .. controls (143.98,418.95) and (143.39,419.55) .. (142.66,419.55) .. controls (141.92,419.55) and (141.33,418.95) .. (141.33,418.22) -- cycle ;
\draw    (150.6,330.14) -- (142.13,421.59) ;
\draw    (230.8,400.74) -- (271.93,421.99) ;
\draw  [dash pattern={on 4.5pt off 4.5pt}]  (282.68,389.74) -- (322.68,411.99) ;
\draw  [dash pattern={on 4.5pt off 4.5pt}]  (230.8,400.74) -- (282.68,389.74) ;
\draw    (271.93,421.99) -- (322.68,411.99) ;
\draw    (280.4,330.54) -- (322.68,411.99) ;
\draw    (280.4,330.54) -- (230.8,400.74) ;
\draw  [color={rgb, 255:red, 245; green, 166; blue, 35 }  ,draw opacity=0.85 ][fill={rgb, 255:red, 245; green, 166; blue, 35 }  ,fill opacity=0.5 ][dash pattern={on 4.5pt off 4.5pt}] (240.5,403.23) -- (273.01,419.25) -- (315.29,410.79) -- (282.78,394.76) -- cycle ;
\draw  [dash pattern={on 4.5pt off 4.5pt}]  (280.4,330.54) -- (282.68,389.74) ;
\draw  [color={rgb, 255:red, 245; green, 166; blue, 35 }  ,draw opacity=1 ][fill={rgb, 255:red, 255; green, 255; blue, 255 }  ,fill opacity=1 ] (239.65,402.74) .. controls (239.65,402.01) and (240.25,401.41) .. (240.98,401.41) .. controls (241.71,401.41) and (242.31,402.01) .. (242.31,402.74) .. controls (242.31,403.48) and (241.71,404.07) .. (240.98,404.07) .. controls (240.25,404.07) and (239.65,403.48) .. (239.65,402.74) -- cycle ;
\draw  [color={rgb, 255:red, 245; green, 166; blue, 35 }  ,draw opacity=1 ][fill={rgb, 255:red, 255; green, 255; blue, 255 }  ,fill opacity=1 ] (280.45,394.43) .. controls (280.45,393.7) and (281.04,393.1) .. (281.78,393.1) .. controls (282.51,393.1) and (283.1,393.7) .. (283.1,394.43) .. controls (283.1,395.16) and (282.51,395.76) .. (281.78,395.76) .. controls (281.04,395.76) and (280.45,395.16) .. (280.45,394.43) -- cycle ;
\draw  [color={rgb, 255:red, 245; green, 166; blue, 35 }  ,draw opacity=1 ][fill={rgb, 255:red, 255; green, 255; blue, 255 }  ,fill opacity=1 ] (313.96,410.79) .. controls (313.96,410.05) and (314.56,409.46) .. (315.29,409.46) .. controls (316.02,409.46) and (316.62,410.05) .. (316.62,410.79) .. controls (316.62,411.52) and (316.02,412.11) .. (315.29,412.11) .. controls (314.56,412.11) and (313.96,411.52) .. (313.96,410.79) -- cycle ;
\draw  [color={rgb, 255:red, 245; green, 166; blue, 35 }  ,draw opacity=1 ][fill={rgb, 255:red, 255; green, 255; blue, 255 }  ,fill opacity=1 ] (271.13,418.62) .. controls (271.13,417.89) and (271.72,417.29) .. (272.46,417.29) .. controls (273.19,417.29) and (273.78,417.89) .. (273.78,418.62) .. controls (273.78,419.35) and (273.19,419.95) .. (272.46,419.95) .. controls (271.72,419.95) and (271.13,419.35) .. (271.13,418.62) -- cycle ;
\draw    (280.4,330.54) -- (271.93,421.99) ;
\draw    (336.92,390.92) .. controls (350.64,385.04) and (353.31,402.21) .. (382.59,396.31) ;
\draw [shift={(384.42,395.92)}, rotate = 167.28] [color={rgb, 255:red, 0; green, 0; blue, 0 }  ][line width=0.75]    (10.93,-3.29) .. controls (6.95,-1.4) and (3.31,-0.3) .. (0,0) .. controls (3.31,0.3) and (6.95,1.4) .. (10.93,3.29)   ;
\draw [color={rgb, 255:red, 245; green, 166; blue, 35 }  ,draw opacity=1 ]   (404,392.84) -- (446.5,412.59) ;
\draw [color={rgb, 255:red, 245; green, 166; blue, 35 }  ,draw opacity=1 ]   (446.5,412.59) -- (485.83,392.17) ;
\draw [color={rgb, 255:red, 245; green, 166; blue, 35 }  ,draw opacity=1 ] [dash pattern={on 4.5pt off 4.5pt}]  (404,392.84) -- (440.97,378.99) ;
\draw [color={rgb, 255:red, 245; green, 166; blue, 35 }  ,draw opacity=1 ] [dash pattern={on 4.5pt off 4.5pt}]  (443.63,378.99) -- (485.83,392.17) ;
\draw  [color={rgb, 255:red, 245; green, 166; blue, 35 }  ,draw opacity=1 ][fill={rgb, 255:red, 245; green, 166; blue, 35 }  ,fill opacity=1 ] (402.67,392.84) .. controls (402.67,392.1) and (403.27,391.51) .. (404,391.51) .. controls (404.73,391.51) and (405.33,392.1) .. (405.33,392.84) .. controls (405.33,393.57) and (404.73,394.16) .. (404,394.16) .. controls (403.27,394.16) and (402.67,393.57) .. (402.67,392.84) -- cycle ;
\draw  [color={rgb, 255:red, 245; green, 166; blue, 35 }  ,draw opacity=1 ][fill={rgb, 255:red, 245; green, 166; blue, 35 }  ,fill opacity=1 ] (444.85,412.59) .. controls (444.85,411.85) and (445.44,411.26) .. (446.17,411.26) .. controls (446.91,411.26) and (447.5,411.85) .. (447.5,412.59) .. controls (447.5,413.32) and (446.91,413.91) .. (446.17,413.91) .. controls (445.44,413.91) and (444.85,413.32) .. (444.85,412.59) -- cycle ;
\draw  [color={rgb, 255:red, 245; green, 166; blue, 35 }  ,draw opacity=1 ][fill={rgb, 255:red, 245; green, 166; blue, 35 }  ,fill opacity=1 ] (484.51,392.17) .. controls (484.51,391.44) and (485.1,390.84) .. (485.83,390.84) .. controls (486.57,390.84) and (487.16,391.44) .. (487.16,392.17) .. controls (487.16,392.9) and (486.57,393.5) .. (485.83,393.5) .. controls (485.1,393.5) and (484.51,392.9) .. (484.51,392.17) -- cycle ;
\draw  [color={rgb, 255:red, 245; green, 166; blue, 35 }  ,draw opacity=1 ][fill={rgb, 255:red, 245; green, 166; blue, 35 }  ,fill opacity=1 ] (440.97,378.99) .. controls (440.97,378.25) and (441.57,377.66) .. (442.3,377.66) .. controls (443.03,377.66) and (443.63,378.25) .. (443.63,378.99) .. controls (443.63,379.72) and (443.03,380.31) .. (442.3,380.31) .. controls (441.57,380.31) and (440.97,379.72) .. (440.97,378.99) -- cycle ;
\draw  [draw opacity=0][fill={rgb, 255:red, 245; green, 166; blue, 35 }  ,fill opacity=0.35 ] (404,391.51) .. controls (403.91,369.74) and (421.58,351.77) .. (443.88,351.19) .. controls (466.48,350.6) and (485.27,368.09) .. (485.85,390.25) .. controls (485.85,390.37) and (485.86,390.49) .. (485.86,390.62) -- (444.93,391.32) -- cycle ; \draw  [color={rgb, 255:red, 245; green, 166; blue, 35 }  ,draw opacity=1 ] (404,391.51) .. controls (403.91,369.74) and (421.58,351.77) .. (443.88,351.19) .. controls (466.48,350.6) and (485.27,368.09) .. (485.85,390.25) .. controls (485.85,390.37) and (485.86,390.49) .. (485.86,390.62) ;  
\draw  [draw opacity=0][fill={rgb, 255:red, 245; green, 166; blue, 35 }  ,fill opacity=0.35 ] (485.82,390.13) .. controls (486.1,411.9) and (468.59,430.02) .. (446.3,430.8) .. controls (423.7,431.58) and (404.76,414.26) .. (403.99,392.1) .. controls (403.98,391.98) and (403.98,391.86) .. (403.98,391.74) -- (444.89,390.67) -- cycle ; \draw  [color={rgb, 255:red, 245; green, 166; blue, 35 }  ,draw opacity=1 ] (485.82,390.13) .. controls (486.1,411.9) and (468.59,430.02) .. (446.3,430.8) .. controls (423.7,431.58) and (404.76,414.26) .. (403.99,392.1) .. controls (403.98,391.98) and (403.98,391.86) .. (403.98,391.74) ;  
\draw    (541.6,374.94) -- (572.77,391.05) ;
\draw  [dash pattern={on 4.5pt off 4.5pt}]  (580.92,366.6) -- (611.23,383.47) ;
\draw  [dash pattern={on 4.5pt off 4.5pt}]  (541.6,374.94) -- (580.92,366.6) ;
\draw    (572.77,391.05) -- (611.23,383.47) ;
\draw    (579.19,321.74) -- (611.23,383.47) ;
\draw    (579.19,321.74) -- (541.6,374.94) ;
\draw  [dash pattern={on 4.5pt off 4.5pt}]  (579.19,321.74) -- (580.92,366.6) ;
\draw    (579.19,321.74) -- (572.77,391.05) ;
\draw    (572.77,391.05) -- (576.01,436.94) ;
\draw    (611.23,383.47) -- (576.01,436.94) ;
\draw    (541.6,374.94) -- (576.01,436.94) ;
\draw  [dash pattern={on 4.5pt off 4.5pt}]  (580.92,366.6) -- (576.01,436.94) ;

\draw (503.73,386.14) node [anchor=north west][inner sep=0.75pt]    {$\sqcup $};
\draw (200.4,378.14) node [anchor=north west][inner sep=0.75pt]    {$\sqcup $};
\draw (170.33,456.41) node [anchor=north west][inner sep=0.75pt]  [font=\large]  {$2\operatorname{Pyr}( \Phi ( P))$};
\draw (339.67,454.74) node [anchor=north west][inner sep=0.75pt]  [font=\large]  {$=$};
\draw (413,455.74) node [anchor=north west][inner sep=0.75pt]  [font=\large]  {$\Phi ( P) \cdot \g{c}\ +\ \operatorname{Bipyr}( \Phi ( P))$};

\end{tikzpicture}

\end{center}
    \caption{Pyramid Operation} \label{fig:pyramid}
\end{figure}
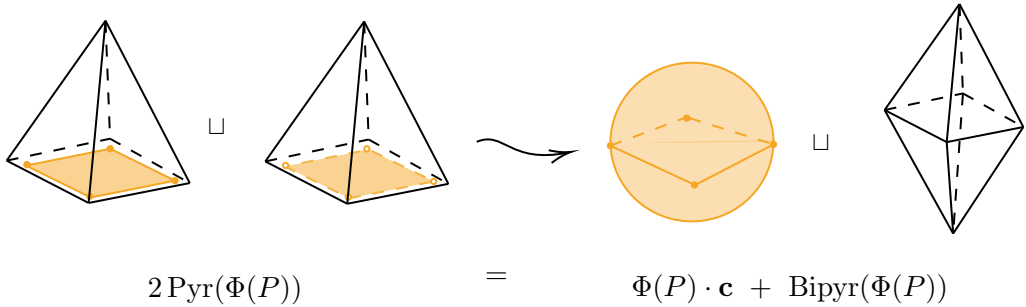

Finally, for the pyramid operation, consider $P_1$ and $P_2$ to be two distinct copies of the polytope $Pyr(P)$. Take a line shelling of $P_1$ (resp. $P_2$) so that the facet corresponding to $P$ is shelled first (resp. last). The contributions of these two facets together will yield $\Phi(P)\cdot \g{c}$, while the remaining contributions will be equal to that of the bipyramid. Equivalently, one can use Corollary 3.14 of Caster--Guyer--Samper~\cite{caster2026} and the fact that the connected sum of two copies of $Pyr(P)$, taken over the facets corresponding to $P$, equals the bipyramid of $P$. Hence one obtains the formula $2\Phi(\text{Pyr}(P))-\Phi(P)\cdot \g{c}=\Phi(\text{Bipyr}(P))$ as desired. We illustrate this operation in Figure~\ref{fig:pyramid}.

\section{Topological Interpretations of $\g{cd}$-monomials}

While the study of the $\g{cd}$-index is often motivated through the lens of
Eulerian posets, we point out that the definition of the $\g{cd}$-index does not require starting from an Eulerian poset. That is, the $\g{cd}$-index can be defined for any abstract collection of chains obeying the Bayer--Billera relations (see \Cref{thm:Bayer-Billera}). 

In this section, we show that every \g{cd}-monomial can be realized as a sub collection of chains of faces of a simple CW complex. To begin, we recall Lemma~1.1 of Stanley~\cite{Stanley1994flag}. Here, we use $\Sigma P$ to denote the \textit{suspension} of $P$, which is combinatorially defined as $P \ast \Delta_1$. Hence, if $P$ is a $d$-sphere, then $\Sigma P$ can be visualized by gluing two $(d+1)$-cells homeomorphically onto $P$.  

\begin{proposition}[Stanley~\cite{Stanley1994flag}]Suppose $P$ and $Q$ are Eulerian. Then, $\Phi_{P \ast Q}(\g{c},\g{d})=\Phi_{P}(\g{c},\g{d})\Phi_{Q}(\g{c},\g{d})$. In particular, 
$\Phi_{\Sigma P}(\g{c},\g{d})=\Phi_P(\g{c},\g{d})\cdot \g{c}$. 
\end{proposition}

Inductively, this recovers the fact that the $\g{cd}$-index of the regular CW  sphere of dimension $d$ with two cells in every (nonnegative) dimension is $\g{c}^{d+1}$. Thus, the question becomes whether one can exhibit topological interpretations of the \g{cd}-monomials that contain at least one $\g{d}$. 

To do this, we consider CW complexes equipped with a cellular $\mathbb{Z}/2$-action. More precisely, we require that that the $\mathbb{Z}/2$-action of $\tau$ on $\Omega$ acts on cells of the same dimension and preserves the poset relation. That is, if $p \preceq q$, then $\tau p \preceq \tau q$. To motivate the relevance of such a $\mathbb{Z}/2$-action, let us consider the circle $S^1$ with $2$ cells in each dimension. Its chains are precisely 
    \[
        \{\emptyset, v_1, v_2, v_1\subseteq e_1, v_1\subseteq e_2, v_2\subseteq e_1,v_2\subseteq e_2\}.
    \]
    Equip this poset with a $\mathbb{Z}/2$-action $\tau$ satisfying $\tau\cdot v_1=v_2$ and $\tau\cdot e_1=e_2$.
    Then, this collection of chains is equipped with a $\mathbb{Z}/2$-action, by extending the action of $\tau$ to each cell of the chain. Consequently, after removing all chains that were fixed (in this case, the empty chain), we obtain the following orbits:
    \[
        \{[v_1], [e_1], [v_1 \subseteq e_1], [v_1 \subseteq e_2]\}.
    \]
     Computing the flag $f$-polynomial obtained from these chains yields $\g{ab} + \g{ba} + 2\g{bb}$. Transforming this into the flag $h$-polynomial, we obtain $\g{ab} + \g{ba}$, which corresponds exactly to the $\g{cd}$-monomial $\g{d}$. 

     We provide an illustration of this process of turning $S^1$ with a \g{cd}-index of $\g{c}^2$ into $(\mathbb{RP}^1,\{\emptyset\})$ with a \g{cd}-index of \g{d} in Figure~\ref{fig:S1toRP1}.
     \begin{figure}[h!]
         \centering
         \begin{tikzcd}
	\bullet & \bullet &&& \bullet \\
	\\
	{\{0\}} & {\{1\}} &&& {\{1\}} \\
	{\{1\}} & {\{0\}} &&& {\{0\}}
	\arrow[""{name=0, anchor=center, inner sep=0}, curve={height=6pt}, no head, from=1-1, to=1-2]
	\arrow[curve={height=-6pt}, no head, from=1-1, to=1-2]
	\arrow["{\mathbb{Z}/2}"', between={0.1}{0.8}, Rightarrow, nfold, from=1-2, to=1-5]
	\arrow[no head, from=1-5, to=1-5, loop, in=55, out=125, distance=10mm]
	\arrow["\Delta", between={0}{0.7}, Rightarrow, nfold, from=1-5, to=3-5]
	\arrow[""{name=1, anchor=center, inner sep=0}, "{\{0,1\}}", no head, from=3-1, to=3-2]
	\arrow["{\{0,1\}}"', no head, from=3-1, to=4-1]
	\arrow["{\{0,1\}}"', no head, from=4-1, to=4-2]
	\arrow[""{name=2, anchor=center, inner sep=0}, "{\{0,1\}}"', no head, from=4-2, to=3-2]
	\arrow["{\{0,1\}}"', curve={height=6pt}, no head, from=4-5, to=3-5]
	\arrow[""{name=3, anchor=center, inner sep=0}, "{\{0,1\}}", curve={height=-6pt}, no head, from=4-5, to=3-5]
	\arrow["\Delta", between={0.2}{0.6}, Rightarrow, nfold, from=0, to=1]
	\arrow["{\mathbb{Z}/2}"', between={0.3}{0.7}, Rightarrow, nfold, from=2, to=3]
\end{tikzcd}
         \caption{Transforming $\g{c}^2$ to $\g{d}$.}
         \label{fig:S1toRP1}
\end{figure}

Next, we can use this construction to consider multiplication on the right by $\g{d}$. Let $\Omega$ be a $d$-sphere. Then, by taking $\Sigma \Sigma \Omega$, one obtains a $(d+2)$-sphere with two $(d+1)$-cells $r_1,r_2$ and two $(d+2)$-cells $\sigma_1,\sigma_2$. This $(d+2)$-sphere $\Sigma \Sigma \Omega$ has a \g{cd}-index of $\Phi(\Omega)\cdot \g{c}^2$. Consider a $\mathbb{Z}/2$-action via $\tau\cdot r_1=r_2$ and $\tau\cdot \sigma_1=\sigma_2$. Quotient all chains by this action and remove all chains contained in $\Omega$. In doing so, we obtain all chains of the form $\mathcal{C}\subseteq [r_1]$, $\mathcal{C}\subseteq [\sigma_1]$, $\mathcal{C}\subseteq [r_1 \subseteq \sigma_1]$ and $\mathcal{C}\subseteq [r_1 \subseteq \sigma_2]$ for every chain $\mathcal{C}$ of $\Omega$. Indeed, this results in the flag $f$-polynomial $\Upsilon_\Omega(\g{a},\g{b})(\g{ab}+\g{ba}+2\g{bb})$, which translates to the \g{cd}-index $\Phi_\Omega(\g{c},\g{d})\cdot \g{d}$ as desired.

Notice that in this operation of multiplying on the right by $\g{d}$, we are always discarding the empty chain, as it is always contained in the chains in the fixed locus of the $\mathbb{Z}/2$-action.

\begin{theorem}
    Every $\g{cd}$-monomial is realizable as the $\g{cd}$-index associated to a sub collection of chains in the face poset of a join of the $0$-sphere and $(\mathbb{RP}^1,\{\emptyset\})$ with one cell in each dimension. 
\end{theorem}

\begin{proof}
    Consider the face poset of the $S^0$ and $\mathbb{RP}^1$ in question where the ``$2$" indicates a $2$ to $1$ gluing map:
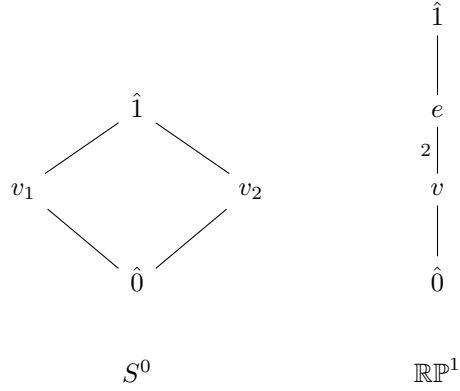
\begin{center}
\begin{tikzcd}
	&&&& {\hat{1}} \\
	& {\hat{1}} &&& e \\
	{v_1} && {v_2} && v \\
	& {\hat{0}} &&& {\hat{0}} \\
	& {S^0} &&& {\mathbb{RP}^1}
	\arrow[no head, from=2-2, to=3-3]
	\arrow[no head, from=2-5, to=1-5]
	\arrow[no head, from=3-1, to=2-2]
	\arrow["2", no head, from=3-5, to=2-5]
	\arrow[no head, from=4-2, to=3-1]
	\arrow[no head, from=4-2, to=3-3]
	\arrow[no head, from=4-5, to=3-5]
\end{tikzcd}
\end{center}
Quickly, one finds the \g{cd}-index of $S^0$ to be \g{c}. Moreover, after removing the empty chain from $\mathbb{RP}^1$, one gets a \g{cd}-index of $\g{d}$. Recall that the contribution of $\g{c}$ in the flag $f$-polynomial is $(\g{a}+2\g{b})$ and the contribution of $\g{d}$ in the flag $f$-polynomial is $(\g{ab}+\g{ba}+2\g{bb})$. Therefore, for a \g{cd}-word \g{w}, consider the poset $\Omega_\g{w}$ which is the poset defined from \g{w} by replacing each $\g{c}$ with the face poset of $S^0$ and each $\g{d}$ with the face poset of $\mathbb{RP}^1$ and multiplication with poset join. Then, in order to obtain the \g{cd}-monomial $\g{w}$, the desired sub collection of chains from $\Omega_\g{w}$ will be the chains that contain at least one element from each copy of $\mathbb{RP}^1$ in this join of posets.
\end{proof}

\begin{example} [A topological interpretation of $\g{d}^2$] To obtain all chains that yield $\g{d}^2$, we consider the following poset: 
\begin{center}
    \begin{tikzcd}
	{\hat{1}} \\
	\sigma \\
	r \\
	e \\
	v \\
	{\hat{0}}
	\arrow[no head, from=2-1, to=1-1]
	\arrow["2", no head, from=3-1, to=2-1]
	\arrow[no head, from=4-1, to=3-1]
	\arrow["2", no head, from=5-1, to=4-1]
	\arrow[no head, from=6-1, to=5-1]
\end{tikzcd}
\end{center}
Then, we consider the sub collection of chains of this poset that contain at least one element of $\{v,e\}$ and at least one element of $\{r,\sigma\}$. Topologically, the weight of $2$ corresponds to a $2$ to $1$ gluing map of cells. Meanwhile, when counting chains, any chain containing both elements will be counted twice e.g., the chain $v \subseteq e \subseteq r \subseteq \sigma$ will be counted with multiplicity four because there are two options for $v\subseteq e$: either $[v_1 \subseteq e_1]$ or $[v_1 \subseteq e_2]$, and two options for $r \subseteq \sigma$: either $[r_1 \subseteq \sigma_1]$ or $[r_1 \subseteq \sigma_2]$. 
\end{example}

\section{Acknowledgments}
The first author thanks Martina Juhnke and Luz Grisales Gomez for helpful conversations. Additionally, the first authors thanks Richard Ehrenborg for the suggestion to try to extend Corollary~\ref{coro:simplexineqs} to cubical shelling components. The authors thank Varun Shah for helpful conversations throughout this project, and thank Margaret Bayer for helpful conversations and feedback on the article. Finally, the authors thank Isabella Novik for thorough feedback during the duration of the work, which included the main idea of the proof of Proposition~3.22. The first author was partially supported by a graduate fellowship from NSF grant DMS-2246399.

\printbibliography
\end{document}